\documentstyle{amsppt}
\TagsOnRight
\NoBlackBoxes
\magnification=\magstep 1
%%%%%%%%%%%%%%%%%%%%%%%%%%%%%%%%%%%%%%%%%%%5
\vsize 23truecm\hsize 16.5truecm
\vbox{\vskip 3truecm}
\parindent6mm
\voffset=-0.5truecm

\topmatter
\title   
On the inverse scattering problem for Jacobi matrices
with the spectrum on an interval, a finite system of
intervals or a Cantor set of positive length 
\endtitle
\rightheadtext{ On the inverse scattering problem for
Jacobi matrices}
\author A. Volberg and  P. Yuditskii 
\endauthor

%\thanks preliminary version
%\endthanks

%%%%%%%%%%%%%
\abstract
Solving inverse scattering problem for a discrete
Sturm--Liouville operator with a rapidly decreasing
potential one gets reflection coefficients $s_\pm$ and
invertible operators
${I+\Cal H_{s_\pm}}$, where $ \Cal H_{s_\pm}$ is the
Hankel operator related to the symbol $s_\pm$. The
Marchenko--Faddeev theorem (in the continuous case)
[6] and
the Guseinov theorem (in the discrete case) [4],
guarantees the uniqueness of solution of the inverse
scattering  problem. In this article we ask the
following natural question --- can one find a precise
condition guaranteeing that the inverse scattering
problem is uniquely solvable and that operators   
${I+\Cal H_{s_\pm}}$
are invertible? Can one claim that uniqueness implies
invertibility or vise versa?

Moreover we are interested here not only in 
the case of 
decreasing potential but also in the case of 
asymptotically
almost periodic potentials. So we merge here two 
mostly developed
cases of inverse problem for Sturm--Liouville
operators:  the
inverse problem with (almost) periodic potential and 
the
inverse problem with the fast decreasing potential.
\endabstract

\endtopmatter

\document

%\rightline{\it Argumentum omni denudatum ornamento}

\head Main results 
\endhead

 %%%%%%%%%%%%%%%%%%%

Asymptotics of polynomials orthogonal on 
a homogeneous set, we described earlier [8],
indicated strongly that there should be a scattering
theory for Jacobi matrix with almost periodic
background like it exists in the classical case of a
constant background. Note that here left and right
asymptotics are not necessary the same almost periodic
coefficient sequences, but they are of the same
spectral class. 
In this work we present all important
ingredients of such theory: reflection/transmission
coefficients, Gelfand--Levitan--Marchenko
transformation operators, a Riemann--Hilbert problem
related to the inverse scattering problem.  At last now
we can say that the reflectionless Jacobi matrices with
homogeneous spectrum are those whose reflection
coefficient is zero.

Moreover, we extend theory in depth and show that 
a reflection coefficient determine uniquely a Jacobi matrix of
the Szeg\"o class and both transformation operators are
invertible  if and only if the spectral density satisfies
matrix $A_2$ condition [10].  

Concerning $A_2$ condition in the inverse scattering 
 we have to mention,
at least as indirect references, the book  
[7, Chapter 2, Sect. 4] and the paper [1].
Generally references to stationary scattering and
inverse scattering problems in connections with spatial
asymptotics can be found in [3], where explicit
expressions of transmission and reflection coefficients
in terms of Weyl functions and phases asymptotic wave
functions was given.
%%%%%%%%%%%%%%%%%%%%%%

Let $J$ be a Jacobi matrix defining a bounded
self--adjoint operator on $l^2(\Bbb Z)$:
$$
J e_n=p_n e_{n-1}+q_n e_{n}+p_{n+1} e_{n+1},
\quad n\in\Bbb Z,\tag 0.1
$$
where $\{e_n\}$ is the standard basis in
$l^2(\Bbb Z)$, $p_n>0$. The resolvent matrix--function
is defined by the relation
$$
R(z)= R(z, J)=\Cal E^*(J-z)^{-1}\Cal E,
\tag 0.2
$$
where $\Cal E:\Bbb C^2 \to l^2(\Bbb Z)$ in such a way
that
$$
\Cal E\bmatrix c_{-1} \\ c_0\endbmatrix=
e_{-1} c_{-1}+e_{0} c_{0}.
$$
This matrix--function possesses an integral
representation
$$
R(z)=\int\frac{d\sigma}{x-z}
\tag 0.3
$$
with a $2\times 2$ matrix--measure having
a compact support on $\Bbb R$. $J$ is unitary
equivalent to the operator multiplication by
independent variable on
$$
L_{d\sigma}^2=
\left\{f=\bmatrix f_{-1}(x) \\
f_0(x)\endbmatrix:\ \int f^*\,d\sigma\,
f<\infty\right\}.
$$
The spectrum of $J$ is called absolutely
continuous if the measure $d\sigma$ is
absolutely
continuous with respect to the Lebesgue measure
on the real axis,
$$
d\sigma(x)=\rho(x)\,dx.\tag 0.4
$$

Let $J_0$ be a Jacobi matrix with constant
coefficients, $p_n=1,\ q_n=0$ (so called
Chebyshev matrix). 
It has the following functional representation,
besides the general one mentioned above.
Note that the resolvent set of $J_0$ is the
domain $\bar\Bbb C\setminus [-2,2]$.
Let $z(\zeta):\Bbb D\to\bar\Bbb C\setminus
[-2,2]$ be a uniformization of this domain,
$z(\zeta)=1/\zeta+\zeta$.
With respect to the standard basis 
$\{t^n\}_{n\in\Bbb Z}$ 
 in
$$
L^2=\{f(t):\ \int_{\Bbb T}|f|^2\,dm\},
$$
the matrix of the operator of multiplication by
$z(t),\ t\in\Bbb T$,
is the
Jacobi matrix $J_0$, since $z(t)
t^n=t^{n-1}+t^{n+1}$.

The famous Bernstein--Szeq\"o theorem implies
the following proposition.

\proclaim{Proposition 0.1} Let $J$ be a Jacobi
matrix whose spectrum is an interval $[-2,2]$.
Assume that the spectrum is absolutely
continuous and the density of the spectral
measure satisfies the condition
$$
\log\det\rho(z(t))\in L^1.
\tag 0.5
$$
Then
$$
p_n\to 1,\ q_n\to 0, \quad n\to\pm\infty.
\tag 0.6
$$
Moreover, there exist  generalized
eigenvectors
$$
\aligned
&p_n e^+(n-1,t)+q_n e^+(n,t)
+p_{n+1} e^+(n+1,t)= z(t) e^+(n,t)\\
&p_n e^-(-n,t)+q_n e^-(-n-1,t)
+p_{n+1} e^-(-n-2,t)= z(t) e^-(-n-1,t)
\endaligned\tag 0.7
$$
such that the following asymptotics hold true
$$
\aligned
s(t)e^{\pm}(n,t)=&s(t) t^n+o(1),\quad
n\to +\infty\\
s(t)e^{\pm}(n,t)=&t^n+s_\mp(t) t^{-n-1}+
o(1),\quad n\to -\infty
\endaligned\tag 0.8
$$
in $L^2$.
\endproclaim

To clarify the meaning of the words
"generalized eigenvectors" we need some
definitions and notation.

The matrix
$$
S(t)=\bmatrix s_-& s\\
               s& s_+
         \endbmatrix(t)\tag 0.9
$$
is called the scattering matrix--function. It is
a unitary--valued matrix--function with the
following
symmetry property:
$$
S^*(\bar t)= S(t),\tag 0.10
$$
and the analytic property:
$$
s(t) 
\text{\ is boundary values of an outer
function.}
\tag 0.11
$$
We still denote by $s(\zeta),\
\zeta\in\Bbb D$,
the values of the function inside the disk,
 and in the sequel, we assume
that $s$ meets the normalization condition
$s(0)>0$.

In fact, these mean that each of the entries
$s_\pm$ (so called reflection coefficient) determines
the matrix
$S(t)$ in unique way. Indeed, since
$$
|s(t)|^2+|s_\pm(t)|^2=1,\tag 0.12
$$
using (0.11), we have
$$
s(\zeta)=e^{\frac 1 2 \int_{\Bbb T}
\frac{t+\zeta}{t-\zeta}\log\{1-|s_\pm(t)|^2\}
\,dm}.
$$
Then, we can solve for $s_\mp$ the relation
$$
\bar s_+ s+\bar s s_-=0.\tag 0.13
$$

With the function $s_\pm$ we associate the
metric
$$\align
||f||^2_{s_\pm}=&
\frac 1 2\left\langle 
\bmatrix 1 &\overline{s_\pm(t)}\\
     s_\pm(t) & 1
\endbmatrix
\bmatrix f(t)\\ \bar t f(\bar t)
\endbmatrix,
\bmatrix f(t)\\ \bar t f(\bar t)
\endbmatrix
\right\rangle\\
=&\langle f(t)+\bar t (s_\pm f)(\bar t),
f(t)\rangle,\quad f\in L^2.
\endalign
$$
Note that the conditions (0.11), (0.12)
guarantee that $||f||_{s_\pm}=0$ implies
$f=0$. We denote by $L^2_{dm,s_\pm}$ or
$L^2_{s_\pm}$ (for shortness) the closer
of $L^2$ with respect to this new metric.

The following relation sets a unitary map
from $L^2_{s_+}$ to $L^2_{s_-}$:
$$
s(t)f^-(t)=\bar t f^+(\bar t) +s_+(t) f^+(t),
$$
moreover, in this case,
$$
||f^+||_{s_+}^2=||f^-||_{s_-}^2=
\frac 1 2
\{||s f^+||^2+||s f^-||^2\},
$$
and the inverse map is of the form
$$
s(t)f^+(t)=\bar t f^-(\bar t) +s_-(t) f^-(t).
$$

We say that a Jacobi matrix $J$ with the spectrum
$[-2,2]$ is of Szeg\"o class if its spectral measure
$d\sigma$ satisfies (0.4), (0.5).

\proclaim{Theorem 0.1} Let $J$ be a Jacobi matrix
 of Szeg\"o class
with the spectrum $E=[-2,2]$. Then
there exists a unique
 unitary--valued matrix--function $S(t)$ of
the form (0.9) possessing the properties
(0.10), (0.11), and a unique pair of Fourier transforms
$$
\Cal F^\pm:l^2(\Bbb Z)\to L^2_{s_\pm},
 \quad
(\Cal F^{\pm} Jf)(t)=z(t)(\Cal F^\pm f)(t),
\tag 0.14
$$
determining each other by the relations
$$
s(t)(\Cal F^\pm f)(t)=
\bar t (\Cal F^\mp f)(\bar t)+
s_\mp(t)(\Cal F^\mp f)(t),
\tag 0.15
$$
and possessing the following analytic properties
$$
s\Cal F^\pm(l^2(\Bbb Z_\pm))\subset H^2,
\tag 0.16
$$
and the asymptotic properties
$$
e^\pm(n,t)=t^n+o(1)\quad\text{in}\ L^2_{s_\pm},
\quad n\to +\infty,
\tag 0.17
$$
where
$$
e^+(n,t)=(\Cal F^+ e_n)(t),\quad
e^-(n,t)=(\Cal F^- e_{-n-1})(t).
$$
(As before, $\{e_n\}$ is the standard basis
in $l^2(\Bbb Z)$).
\endproclaim

\remark
{ Remark 0.1} Show that (0.17) is equivalent
to (0.8). Due to 
$$
\bmatrix 1 & \bar s_\pm\\
         s_\pm& 1
\endbmatrix=
\bmatrix |s|^2 & 0\\
         0 & 0
\endbmatrix+
\bmatrix \bar s_\pm\\  1
\endbmatrix
\bmatrix 
         s_\pm& 1
\endbmatrix,
\tag 0.18
$$
(0.17) is equivalent to ($n\to +\infty$)
$$\aligned
&s(t) e^{\pm}(n,t)=s(t) t^n+o(1)\quad \text{in}\  L^2,\\
&s_\pm(t) e^{\pm}(n,t)
+\bar t e^{\pm}(n,\bar t)=s_\pm(t) t^n+
\bar t^{n+1}+ o(1)\quad \text{in}\  L^2.
\endaligned
$$
Using (0.15), we rewrite the second relation into the
form
$$
s(t) e^{\mp}(-n-1,t)= t^{-n-1}+
s_\pm(t) t^n+
o(1)\quad \text{in}\  L^2.
$$
Substituting $n:=-n-1$, we get the second relation
of (0.8).

\endremark

\bigskip

 A fundamental question is how to recover the Jacobi
matrix from the scattering matrix, in fact, from the
reflection coefficient $s_+$ (or $s_-$)? When can this
be done? Do we have a uniqueness theorem?

We show that for an arbitrary function $s_+(t)$ 
satisfying
$$
\overline{s_+(\bar t)}=s_+(t)\quad\text{and}\quad
\log\{1-|s_+(t)|^2\}\in L^1,
\tag 0.19
$$
there exists a Jacobi matrix $J$   of
Szeg\"o class  such  that $s_+(t)$ is its
reflection coefficient. But we can construct
a matrix with this property, at least,
in two different ways.

First, consider the space
$$
H^2_{s_+}=\text{clos}_{L^2_{s_+}} H^2,
$$
and introduce the Hankel operator $\Cal H_{s_+}: H^2\to
H^2$,
$$
\Cal H_{s_+} f= P_+\bar t(s_+ f)(\bar t),
\quad f\in H^2,
$$
where $P_+$ is the Riesz projection from $L^2$
onto $H^2$. This operator determines the metric
in $H^2_{s_+}$:
$$\aligned
||f||^2_{s_+}=&\langle f(t)+\bar t(s_+ f)(\bar t),
f(t)\rangle\\
=&\langle (I+\Cal H_{s_+}) f, f\rangle,\quad
\forall f\in H^2.
\endaligned
$$

\proclaim{Lemma 0.1} Under the assumptions (0.19),
the space $ H_{s_+}^2$ is  a space of holomorphic
functions with a reproducing kernel. Moreover,
$s f\in H^2$ for any $f\in H^2_{s_+}$, and
the reproducing vector $k_{s_+}$:
$$
\langle f, k_{s_+}\rangle=f(0),\quad\forall f\in 
 H_{s_+}^2,
$$
is of the form
$$
 k_{s_+}=(I+\Cal H_{s_+} )^{[-1]} \text{\bf 1}:=
\lim_{\epsilon\to 0^+}(\epsilon+I+\Cal H_{s_+} )^{-1}
\text{\bf 1}
\quad\text{in}\ L^2_{s_+}.\tag 0.20
$$
\endproclaim

Put $K_{s_+}(t)=k_{s_+}(t)/\sqrt{k_{s_+}(0)}$.

\proclaim{Theorem 0.2} Let $s_+(t)$ satisfy (0.19).
Then the system of functions 
$\{t^n K_{s_+t^{2n}}(t)\}_{n\in\Bbb Z}$ forms an
orthonormal basis in $L^2_{s_+}$. With respect to
this basis, operator multiplication by $z(t)$ is
a Jacobi matrix $J$ of  Szeg\"o class. Moreover,
the initial function $s_+(t)$ is the reflection
coefficient of the scattering matrix--function
$S(t)$, associated to $J$ by Theorem 0.1, and
$$
e^+(n,t)= t^n K_{s_+t^{2n}}(t).
$$
\endproclaim

From the other hand, the system of functions
$\{t^n K_{s_-t^{2n}}(t)\}_{n\in\Bbb Z}$ forms an
orthonormal basis in $L^2_{s_-}$, and we are able
to define a Jacobi matrix $\tilde J$ by the relation
$$
z(t)\tilde e^+(n,t)=
\tilde p_n \tilde e^+(n-1,t)+\tilde q_n 
\tilde e^+(n,t)
+\tilde p_{n+1} \tilde e^+(n+1,t),
$$ 
where $\{\tilde e^+(n,t)\}$ is the dual system
to the system $\{t^n K_{s_-t^{2n}}(t)\}$
(see (0.15)), i.e.:
$$
s(t)\tilde e^+(-n-1,t)=\bar t^{n+1} K_{s_-t^{2n}}
(\bar t)+ s_-(t) t^n K_{s_-t^{2n}}(t).
$$
Even the invertibility condition for the operators
$(I+\Cal H_{s_\pm})$ does not guarantee that operators
$J$ and $\tilde J$ are the same. But if $J=\tilde J$,
then the uniqueness theorem takes place.

\proclaim{Theorem 0.3} Let $s_+$ satisfy (0.19). 
Then the reflection coefficient $s_+$ determines a 
Jacobi matrix $J$ of Szeg\"o class in a unique way
if and only if
the following relations take place
$$
s(0)K_{s_\pm}(0) K_{s_\mp t^{-2}}(0)=1.\tag 0.21
$$
\endproclaim

\proclaim{Corollary 0.1} Let $J$ be a Jacobi matrix
of Szeg\"o class with the spectrum $[-2,2]$ and let
$\rho$ be the density of its spectral measure. If
$$
\int_{-2}^2\rho^{-1}(x)\,dx<\infty,
$$
then there is no other Jacobi matrix of Szeg\"o
class with the same scattering matrix--function $S(t)$.
\endproclaim

It is important to know, when the operators 
$(I+\Cal H_{s_\pm})$,  
playing a central role in 
the inverse scattering problem, are invertible in the  
proper sense of words.

\proclaim{Theorem 0.4}
 Let $J$ be a Jacobi matrix of
Szeq\"o class  with the spectrum $[-2,2]$. Let 
$\rho$ be
the density of its spectral measure and let $s_+$ be the
reflection coefficient of its scattering
matrixÑ-function. Then the following statements are
equivalent.

\item
{1.} The spectral density $\rho$ satisfies condition
$A_2$. 
\item
{2.} The reflection coefficient $s_+$ determines  a 
Jacobi
 matrix of Szeq\"o class uniquely and both operators 
$(I+\Cal H_{s_\pm})$ are invertible. 
\endproclaim

\bigskip

To extend these results to the case when a spectrum
$E$ is a finite system of intervals or a Cantor set
of positive measure, we  need only to introduce a
counterpart of Hardy space.

Let $z(\zeta):\Bbb D\to \Omega$ be a uniformization of
the domain $\Omega=\bar\Bbb C\setminus E$.
Thus there exists  
a discrete subgroup $\Gamma$ of the group $SU(1,1)$ 
consisting of elements of the form
$$
\gamma=\bmatrix \gamma_{11}&\gamma_{12}\\ \gamma_{21}
&\gamma_{22}
\endbmatrix,\ 
\gamma_{11}=\overline{\gamma_{22}},\ \gamma_{12}=
\overline{\gamma_{21}},
\ \det\gamma=1,
$$
such that  $z(\zeta)$ 
is automorphic with respect to $\Gamma$, i.e.,
$z(\gamma(\zeta))=z(\zeta),\ \forall \gamma\in\Gamma$,
and any two preimages of $z_0\in\Omega$ are 
$\Gamma$--equivalent, i.e.,
$$
z(\zeta_1)=z(\zeta_2)\ \Rightarrow\ \exists\gamma
\in\Gamma:\ 
\zeta_1=\gamma(\zeta_2).
$$
We normalize $z(\zeta)$ by the conditions $z(0)=\infty$,
$(\zeta z)(0)>0$.

A character of $\Gamma$ is a complex--valued function
$\alpha:\Gamma\to\Bbb T$, satisfying
$$
\alpha(\gamma_1\gamma_2)=\alpha(\gamma_1)
\alpha(\gamma_2),\quad\gamma_1,\gamma_2\in\Gamma.
$$
The characters form an Abelian compact group denoted
by $\Gamma^*$.

For a given character $\alpha\in\Gamma^*$, as usual let
us define 
$$
H^\infty(\Gamma,\alpha)= \{f\in H^\infty:
\ f(\gamma(\zeta))=\alpha(\gamma)f(\zeta),\ \forall
\gamma\in\Gamma\}.
$$
Generally, a group
$\Gamma$ is said to be of {\it Widom type} if for any
$\alpha\in\Gamma^*$ the space $H^\infty(\Gamma,\alpha)$
is not trivial (contains a non--constant function).

A group of Widom type  acts
dissipatively on
$\Bbb T$ with respect to $dm$, that is there exists a 
measurable (fundamental)
set $\Bbb E$, which does not contain any two 
$\Gamma$--equivalent points,
and the union $\cup_{\gamma\in\Gamma}\gamma(\Bbb E)$ 
is a set of 
full measure. We can choose $\Bbb E$ possessing the
symmetry property: $t\in \Bbb E\Rightarrow\bar t\in \Bbb
E$. For the space of square summable functions 
on $\Bbb E$ (with respect
to the Lebesgue measure), we use
the notation $L^2_{dm|\Bbb E}$.

Let $f$ be an analytic function in 
$\Bbb D$, $\gamma\in\Gamma$ and $k\in \Bbb N$.
Then we put
$$
f\vert[\gamma]_k=\frac{f(\gamma(\zeta))}
{(\gamma_{21}\zeta+\gamma_{22})^k}.
$$
Notice that $f\vert[\gamma]_2=
f\ \forall\gamma\in \Gamma$, means 
that the form
$f(\zeta)d\zeta$ is invariant with respect to 
the substitutions
$\zeta\to\gamma(\zeta)$ ($f(\zeta)d\zeta$ is an Abelian
integral on $\Bbb D/\Gamma$). Analogically, 
$f\vert[\gamma]=
\alpha(\gamma)f\ \forall\gamma\in \Gamma$, means 
that the form
$|f(\zeta)|^2\,|d\zeta|$ is invariant with respect to 
these substitutions.

We recall, that a function $f(\zeta)$ is of Smirnov
class, if it can be represented as a ratio of two
functions from $H^\infty$ with an outer denominator.

\proclaim{Definition} 
Let $\Gamma$ be a group of Widom type. 
The space $A^{2}_1(\Gamma,\alpha)$
($A^{1}_2(\Gamma,\alpha)$) is formed by functions $f$,
which are analytic on $\Bbb D$ and satisfy the
following three conditions
$$
\align
1)& f \ \text{is of Smirnov class}\\
2)& f\vert[\gamma]=\alpha(\gamma) f\ \ \ 
(f\vert[\gamma]_2=\alpha(\gamma) f)
\quad 
\forall\gamma\in \Gamma\\
3)& 
\int_{\Bbb E}\vert f\vert^2\,dm<\infty\  \ \ 
(\int_{\Bbb E}\vert f\vert\,dm<\infty).
\endalign
$$
\endproclaim

$A^2_1(\Gamma,\alpha)$
is a Hilbert space with the reproducing kernel
$k^\alpha(\zeta,\zeta_0)$, moreover
$$
0<\inf_{\alpha\in\Gamma^*} k^\alpha(\zeta_0,\zeta_0)\le
\sup_{\alpha\in\Gamma^*} k^\alpha(\zeta_0,\zeta_0)
<\infty.\tag W
$$
Put
$$
k^\alpha(\zeta)=k^\alpha(\zeta,0)\quad\text{and}
\quad
K^\alpha(\zeta)=\frac{k^\alpha(\zeta)}
{\sqrt{k^\alpha(0)}}.
$$

We need one more special function.
 The Blaschke product
$$
b(\zeta)=\zeta\prod_{\gamma\in\Gamma, \gamma\not= 1_2} 
\frac{\gamma(0)-\zeta}{1-\overline{\gamma(0)}\zeta}
 \frac{\vert\gamma(0)\vert}{\gamma(0)}
$$
is called the {\it Green's function} of $\Gamma$ with 
respect to the origin.
It is character--automorphic function, i.e., there
exists $\mu\in\Gamma^*$ such that $b\in
H^\infty(\Gamma,\mu)$.
 Note, if $G(z)=G(z,\infty)$ denotes 
the Green's function of the domain
$\Omega$, then
$$
G(z(\zeta))=-\log\vert b(\zeta)\vert.
$$
 
\proclaim{Theorem [5]}
Let $\Gamma$ be a group of Widom type. 
The following statements are equivalent:

\noindent 1) The function $K^\alpha(0)$ is continuous on
$\Gamma^*$.

\noindent 2)  
$
\sup\{\vert f(0)\vert:
f\in
H^\infty(\Gamma,\alpha),\ \Vert f\Vert\le 1\}\to 1
$,\quad
 $\alpha\to 1_{\Gamma^*}$. 

\noindent 3)
The {\it Direct Cauchy Theorem} holds:
$$
\int_{\Bbb E}\frac{f}{b}(t)\,\frac{d t}{2\pi
i} =\frac{f}{b'}(0),
\quad\forall f\in A^1_2(\Gamma,\mu).\tag DCT
$$

\noindent 4) Let 
$\overline{t  A^2_1(\Gamma,\alpha^{-1})}=
 \{g=\overline{t f}:\ f\in A^2_1(\Gamma,\alpha^{-1})\}$. Then
$$
L^2_{dm\vert\Bbb E} =
\overline{t  A^2_1(\Gamma,\alpha^{-1})}\oplus
 A^2_1(\Gamma,\alpha)\quad \forall\alpha\in\Gamma^*.
 $$
 
 \noindent 5) Every invariant subspace  
$M\subset A^2_1(\Gamma,\alpha)$
(i.e. $\phi M\subset M\ \forall \phi\in
H^\infty(\Gamma)$) is of the form
$$ M=\Delta A^2_1(\Gamma,\beta^{-1}\alpha)$$ for some
character--automorphic inner function $\Delta\in
H^\infty(\beta)$.

\endproclaim

\proclaim{Definition [2]}
A measurable set $E$ is homogeneous if there is an 
$\eta>0$ such that
$$
\vert(x-\delta,x+\delta)\cap E\vert\ge 
\eta\delta\quad \text{for all}\
0<\delta<1 \quad\text {and all}\ 
x\in E.\tag C
$$
\endproclaim

A Cantor set of  positive length is an example of
a homogeneous set. Let $E$ be a homogeneous set,
then the domain $\Omega=\bar C\setminus E$ (respectively
the group $\Gamma$) is of Widom type and the Direct
Cauchy Theorem holds.

\bigskip
 Recall, that a sequence
of real numbers $\{p_n\}\in l^\infty(\Bbb Z)$ is called
uniformly almost periodic if the set of sequences
$\{\{p_{n+l}\},\ l\in\Bbb Z\}$ is a precompact in
$l^\infty(\Bbb Z)$. The general way to produce a
sequence of this type looks as follows: let
$\Cal G$ be a compact Abelian group, and let
$f(g)$ be a continuous function on $\Cal G$,
then
$$
p_n:=f(g_0+n g_1),\quad g_0, g_1\in\Cal G,
$$
is an almost periodic sequence. A Jacobi matrix is
almost periodic  if the coefficient sequences are
almost periodic.  We denote by $J(E)$ the class of
almost periodic Jacobi matrices with absolutely
continuous homogeneous spectrum
$E$. In what follows the class $J(E)$ plays a role of
Chebyshev matrix. In fact,
if $E=[-2,2]$ then $J(E)=\{J_0\}$.

\proclaim{Theorem [9]} Let $E$ be a homogeneous set.
Let $z:\Bbb D\to\bar \Bbb C\setminus E$ be a
uniformizing mapping.  
Then the
systems of functions 
$\{b^n K^{\alpha\mu^{-n}}\}_{n\in\Bbb Z_+}$ 
and $\{b^n K^{\alpha\mu^{-n}}\}_{n\in\Bbb Z}$ 
form an
orthonormal basis in $A^2_1(\Gamma,\alpha)$ and
 in $L^2_{dm\vert\Bbb E} $, respectively,  
for any $\alpha\in\Gamma^*$.
With respect to this basis, the operator multiplication 
by $z(t)$ is a
three--diagonal almost periodic Jacobi matrix
$J(\alpha)$. Moreover,
$$
 J(E)=\{J(\alpha):\ \alpha\in\Gamma^*\},
$$
and $J(\alpha)$ is a continuous function on $\Gamma^*$. 
\endproclaim
%%%%%%%%%%%%%%%%%%%%%%%%%

We say that a Jacobi matrix $J$ with the spectrum
$E$ is of Szeg\"o class if its spectral measure
is absolutely continuous, 
$d\sigma(x)=\rho(x)\,dx$, and $\rho(z(t))$ satisfies
(0.5).

\proclaim{Theorem 0.5} Let $J$ be a Jacobi matrix
 of Szeg\"o class
with a homogeneous spectrum $E$. Then
there exists a unique
 unitary--valued matrix--function $S(t)$ of
the form (0.9) possessing the properties
(0.10), (0.11), and a unique pair of Fourier transforms
$$
\Cal F^\pm:l^2(\Bbb Z)\to L^2_{dm|\Bbb E, s_\pm},
 \quad
(\Cal F^{\pm} Jf)(t)=z(t)(\Cal F^\pm f)(t),
\tag 0.22
$$
determining each other by the relations
$$
s(t)(\Cal F^\pm f)(t)=
\bar t (\Cal F^\mp f)(\bar t)+
s_\mp(t)(\Cal F^\mp f)(t),
\tag 0.23
$$
and possessing the following analytic properties
$$
s\Cal F^\pm(l^2(\Bbb Z_\pm))\subset
A_1^2(\Gamma,\alpha^{-1}_\mp),
\tag 0.24
$$
and the asymptotic properties
$$
e^\pm(n,t)=b^n(t)K^{\alpha_\pm\mu^{-n}}(t)+
o(1)\quad\text{in}\
L^2_{dm|\Bbb E, s_\pm},
\quad n\to +\infty,
\tag 0.25
$$
where
$$
e^+(n,t)=(\Cal F^+ e_n)(t),\quad
e^-(n,t)=(\Cal F^- e_{-n-1})(t),
$$
and $L^2_{dm|\Bbb E, s_\pm}$ is the closer
of the functions from $L^2_{dm|\Bbb E}$ with respect
to the metric
$$
||f||^2_{s_\pm}=
\frac 1 2\left\langle 
\bmatrix 1 &\overline{s_\pm(t)}\\
     s_\pm(t) & 1
\endbmatrix
\bmatrix f(t)\\ \bar t f(\bar t)
\endbmatrix,
\bmatrix f(t)\\ \bar t f(\bar t)
\endbmatrix
\right\rangle,
\quad f\in L^2_{dm|\Bbb E}.
$$
 
 \endproclaim

Theorems 0.2--0.4 also have their closely parallel
counterparts in the case when the spectrum is a
homogeneous set.

%%%%%%%%%%%%%%%%%%
%We finish this manuscript with a remark 
%on a connection between this new type 
%inverse scattering
%problem and the RiemannÑ-Hilbert problem.
%%%%%%%%%%%%%%%%%%%

 %%%%%%%%%%%%%%%%%%%%%%%%%%%

\head In the model space
\endhead

Let $E$ be a homogeneous set. Let 
$z(\zeta):\Bbb D/\Gamma\sim \bar\Bbb C\setminus E$ be a
uniformization and
$b(\zeta)$ be the Green's function.  Throughout the 
paper we assume that $(bz)(0)=1$. Let
$\Bbb E\subset \Bbb T$ be a symmetric fundamental set
($t\in\Bbb E\Rightarrow\bar t\in \Bbb E$). 

With a function
$s_+(t)\in L^\infty_{dm|\Bbb E}$ such that
$$
\overline{s_+(\bar t)}=s_+( t)\quad\text{and}\quad
1-|s_+(t)|^2>0\ \ \text{a.e. on}\ \Bbb E,
\tag 1.1
$$
we associate the metric
$$\align
||f||^2_{s_+}=&
\frac 1 2\left\langle 
\bmatrix 1 &\overline{s_+(t)}\\
     s_+(t) & 1
\endbmatrix
\bmatrix f(t)\\ \bar t f(\bar t)
\endbmatrix,
\bmatrix f(t)\\ \bar t f(\bar t)
\endbmatrix
\right\rangle\\
=&\langle f(t)+\bar t (s_+ f)(\bar t),
f(t)\rangle,\quad f\in L^2_{dm|\Bbb E}.
\endalign
$$
Condition (1.1) guarantee that $||f||_{s_+}=0$
implies $f=0$. We denote by $L^2_{dm|\Bbb E,s_+}$
or $L^2_{s_+}$ (for shortness) the closer of
$L^2_{dm|\Bbb E}$ with respect to this metric.

\proclaim{Lemma 1.1} The operator multiplication
by $z(t)$ in $L^2_{s_+}$ is unitary equivalent to
the operator multiplication by $z(t)$ in
$L^2_{dm|\Bbb E}$
\endproclaim

\demo{Proof} Let us put
$$
\bmatrix g(t)\\ \bar t g(\bar t)\endbmatrix=
\bmatrix 1 &\overline{s_+(t)}\\
     s_+(t) & 1
\endbmatrix^{1/2}
\bmatrix f(t)\\ \bar t f(\bar t)
\endbmatrix,\quad
f\in L^2_{dm|\Bbb E}.
$$
In this case $||f||_{s_+}=||g||$. The system of 
identities
$$
\align
\bmatrix 1 &\overline{s_+(t)}\\
     s_+(t) & 1
\endbmatrix^{1/2}
\bmatrix (zf)(t)\\ \bar t (zf)(\bar t)
\endbmatrix=&
\bmatrix 1 &\overline{s_+(t)}\\
     s_+(t) & 1
\endbmatrix^{1/2} z(t)
\bmatrix f(t)\\ \bar t f(\bar t)
\endbmatrix\\
=&z(t)\bmatrix g(t)\\ \bar t g(\bar t)\endbmatrix=
\bmatrix (zg)(t)\\ \bar t (zg)(\bar t)\endbmatrix
\endalign
$$
finishes the proof.\qed
\enddemo

Let $\alpha_+\in\Gamma^*$. In what further, we assume 
that
$s_+\in L^\infty(\Gamma,\alpha_+^{-2})$ and
$$
\log(1-|s_+(t)|^2)\in L^1.\tag 1.2
$$
We define an outer function $s$, $s(0)>0$, by the 
relation
$$
|s(t)|^2=1-|s_+(t)|^2,\quad t\in\Bbb T.
$$
It is a character--automorphic function
such that $\overline{s(\bar t)}=s(t)$. It is convenient
to denote its character by $\alpha_+^{-1}\alpha_-^{-1}$,
i.e., $s\in H^\infty(\Gamma,\alpha_+^{-1}\alpha_-^{-1})$.

Let us discuss some properties of the space
$$
H^2_{s_+}(\alpha_+)
:=\text{clos}_{L^2_{s_+}} A_1^2(\Gamma,
\alpha_+).
$$
First of all, we define "a Hankel operator"
$\Cal H_{s_+}:A_1^2(\Gamma,\alpha_+)\to 
A_1^2(\Gamma,\alpha_+)$,
$$
\Cal H_{s_+}f=P_{A_1^2(\Gamma,\alpha_+)}
\bar t(s_+ f)(\bar t).
$$
Note, that this operator, indeed, does not
depend 
on "an analytical part" of its symbol, more precisely,
$$
\Cal H_{(s_++\epsilon)}=\Cal H_{s_+},
\quad\forall\epsilon\in H^\infty(\Gamma,\alpha^{-2}_+).
$$
Besides, in the classical case $E=[-2,2]$, 
$\Gamma=\{1_2\}$,
$\Bbb E=\Bbb T$, with a function
$$
s_+(t)=\sum_{n\in\Bbb Z} a_n t^n
$$
is associated the operator $\Cal H_{s_+}:H^2\to H^2$
having the representation
$$
\Cal H_{s_+}=\bmatrix
a_{-1}&a_{-2}&a_{-3}&\dots\\
a_{-2}&a_{-3}&\dots&\\
a_{-3}&\dots& &\\
\dots& & & 
\endbmatrix
$$
with respect to the standard basis $\{t^n\}_{n\in 
\Bbb Z_+}$
in $H^2$. 

The operator $\Cal H_{s_+}$ determines the metric in
$H^2_{s_+}(\alpha_+)$:
$$\align
||f||^2_{s_+}=&\langle f(t)+\bar t (s_+ f)(\bar t),
f(t)\rangle\\
=&\langle (I+ \Cal H_{s_+})f,
f \rangle,\quad f\in A_1^2(\Gamma,\alpha_+).
\endalign
$$

\proclaim{Lemma 1.2} Under the assumptions (1.2),
the space $ H_{s_+}^2(\alpha_+)$ is  a space of 
holomorphic
functions with a reproducing kernel. Moreover,
$s f\in A_1^2(\Gamma,\alpha_-^{-1})$ for any 
$f\in H^2_{s_+}(\alpha_+)$,
and the reproducing vector $k_{s_+}^{\alpha_+}$:
$$
\langle f, k_{s_+}^{\alpha_+}\rangle=f(0),\quad
\forall f\in 
 H_{s_+}^2(\alpha_+),
$$
is of the form
$$
 k_{s_+}^{\alpha_+}=(I+\Cal H_{s_+} )^{[-1]} 
k^{\alpha_+}:=
\lim_{\epsilon\to 0^+}(\epsilon+I+\Cal H_{s_+} )^{-1}
k^{\alpha_+}
\quad\text{in}\ L^2_{s_+}.\tag 1.3
$$
 
\endproclaim

\demo{Proof} From the inequality
$$
\bmatrix 1 &\overline{s_+(t)}\\
     s_+(t) & 1
\endbmatrix
-\bmatrix |s(t)|^2 &0\\
     0 & 0
\endbmatrix=
 \bmatrix |s_+(t)|^2 &\overline{s_+(t)}\\
     s_+(t) & 1
\endbmatrix\ge 0,\tag 1.4
$$
it follows that
$$
||sf||^2\le 2||f||^2_{s_+}\quad\forall f\in L^2_{s_+}.
$$
Thus, if a sequence $\{f_n\}$, $f_n\in
A_1^2(\Gamma,\alpha_+)$, converges in $H^2_{s_+}
(\alpha_+)$,
then the sequence $\{s f_n\}$ converges in
$A_1^2(\Gamma,\alpha_-^{-1})$. In the same way we have
boundedness of the functional
$f\to f(0)$,
$$
|f(0)|^2\le\frac 1{|s(0)|^2}|(sf)(0)|^2\le
\frac 2{|s(0)|^2}||f||^2_{s_+}k^{\alpha_-^{-1}}(0).
$$

Let us prove (1.3). Let $\epsilon>0$, then for 
the norm of
the difference we have an estimate
$$\align
||k_{s_+}^{\alpha_+}-&
(\epsilon+I+\Cal H_{s_+} )^{-1}k^{\alpha_+}||_{s_+}^2=
k_{s_+}^{\alpha_+}(0)-2\{(\epsilon+I+\Cal H_{s_+}
)^{-1}k^{\alpha_+}\}(0)\\
+&\langle(I+\Cal H_{s_+})
(\epsilon+I+\Cal H_{s_+}
)^{-1}k^{\alpha_+},
(\epsilon+I+\Cal H_{s_+}
)^{-1}k^{\alpha_+}\rangle\\
\le&k_{s_+}^{\alpha_+}(0)-\{(\epsilon+I+\Cal H_{s_+}
)^{-1}k^{\alpha_+}\}(0).\tag 1.5
\endalign
$$
Therefore,
$$
\{(\epsilon+I+\Cal H_{s_+}
)^{-1}k^{\alpha_+}\}(0)\le k_{s_+}^{\alpha_+}(0).
\tag 1.6
$$
Besides, (1.5) implies that (1.3) follows from the
relation
$$
\lim_{\epsilon\to 0}\{(\epsilon+I+\Cal H_{s_+}
)^{-1}k^{\alpha_+}\}(0)= k_{s_+}^{\alpha_+}(0).
\tag 1.7
$$

Let us prove (1.7). Since the function
$$
\{(\epsilon+I+\Cal H_{s_+}
)^{-1}k^{\alpha_+}\}(0)=
\langle (\epsilon+I+\Cal H_{s_+}
)^{-1}k^{\alpha_+}, k^{\alpha_+}\rangle
$$
decreases with $\epsilon$ and it is bounded by (1.6),
there exists a limit
$$
\lim_{\epsilon\to 0}\{(\epsilon+I+\Cal H_{s_+}
)^{-1}k^{\alpha_+}\}(0)\le k_{s_+}^{\alpha_+}(0).
\tag 1.8
$$

From the other hand, for any 
$f\in A^2_1(\Gamma,\alpha_+)$
and $\epsilon>0$ the following inequalities hold
$$\align
|f(0)|^2\le&
\langle (\epsilon+I+\Cal H_{s_+})^{-1}k^{\alpha_+},
k^{\alpha_+}\rangle
\langle (\epsilon+I+\Cal H_{s_+})f,f\rangle\\
\le&\{\lim_{\epsilon\to 0}
\langle (\epsilon+I+\Cal H_{s_+})^{-1}k^{\alpha_+},
k^{\alpha_+}\rangle\}
\langle (\epsilon+I+\Cal H_{s_+})f,f\rangle,
\endalign
$$
that is
$$
|f(0)|^2\le
\{\lim_{\epsilon\to 0}
\langle (\epsilon+I+\Cal H_{s_+})^{-1}k^{\alpha_+},
k^{\alpha_+}\rangle\}
||f||^2_{s_+}.
$$
Putting $f=k^{\alpha_+}_{s_+}$, we have
$$
k^{\alpha_+}_{s_+}(0)\le
\lim_{\epsilon\to 0}
\langle (\epsilon+I+\Cal H_{s_+})^{-1}k^{\alpha_+},
k^{\alpha_+}\rangle.
$$
Comparing this inequality with (1.8), we get (1.7), 
thus
(1.3) is proved.\qed
\enddemo

\medskip
We define $s_-\in L^\infty(\Gamma,\alpha^{-2}_-)$ by 
$$
s_-(t)=-\overline{s_+(t)} s(t)/\overline{s(t)}.
$$
In this case 
$$
S(t)=\bmatrix s_-& s\\
s& s_+\endbmatrix(t)
$$
is a unitary--valued matrix function possessing 
properties
(0.10), (0.11).

\proclaim{ Lemma 1.3}
The following relation sets a unitary map
from $L^2_{s_+}$ to $L^2_{s_-}$:
$$
s(t)f^-(t)=\bar t f^+(\bar t) +s_+(t) f^+(t).
$$
Moreover, in this case,
$$
||f^+||_{s_+}^2=||f^-||_{s_-}^2=
\frac 1 2
\{||s f^+||^2+||s f^-||^2\},
$$
and the inverse map is of the form
$$
s(t)f^+(t)=\bar t f^-(\bar t) +s_-(t) f^-(t).
$$
\endproclaim
\demo{Proof} These follow from the identities
$$
\bmatrix 1& \bar s_+\\
s_+& 1\endbmatrix=
\bmatrix \bar s_+& 1\\
1& s_+\endbmatrix
\bmatrix 1/\bar s & 0\\
0& 1/s\endbmatrix
\bmatrix 1& \bar s_-\\
s_-& 1\endbmatrix
\bmatrix 1/s& 0\\
0& 1/\bar s\endbmatrix
\bmatrix s_+& 1\\
1& \bar s_+\endbmatrix
$$
and (0.18).\qed
 
\enddemo

\proclaim{Lemma 1.4} Let $K^{\alpha_+}_{s_+}(t)=
k^{\alpha_+}_{s_+}(t)/ \sqrt{k^{\alpha_+}_{s_+}(0)}$.
The system of functions 
$\{b^n(t)K^{\alpha_+\mu^{-n}}_{s_+b^{2n}}(t)\}$ forms an
orthonormal basis in 
$H^2_{s_+}(\alpha_+)$ when
$\{n\in\Bbb Z_+\}$ and in $L^2_{s_+}$ when
$\{n\in\Bbb Z\}$. With respect to this basis the 
operator
 multiplication by $z(t)$ is a Jacobi matrix.
\endproclaim

\demo{Proof} First, we note that
$$
\{f:\ f\in  H^2_{s_+}(\alpha_+),\ f(0)=0\}=
\{f=b\tilde f:\ \tilde f\in H^2_{s_+b^2}
(\alpha_+\mu^{-1})\}.
$$
Therefore,
$$
 H^2_{s_+}(\alpha_+)=\{K^{\alpha_+}_{s_+}(t)\}
\oplus b H^2_{s_+b^2}(\alpha_+\mu^{-1}).
$$
Iterating this relation, we get that
$\{b^n(t)K^{\alpha_+\mu^{-n}}_{s_+b^{2n}}(t)\}_{n\in\Bbb
Z_+}$ is an orthonormal basis in $H^2_{s_+}(\alpha_+)$,
since $\cap_{n\in\Bbb Z_+}b^n 
H^2_{s_+b^{2n}}(\alpha_+\mu^{-n})=\{0\}$.

Then, we note that an arbitrary function 
$f\in L^2_{s_+}$
can be approximated with the given accuracy by a 
function
$f_1$ from $L^2_{dm|\Bbb E}$. This function, in its 
turn,
can be approximated by a function 
$f_2\in b^n A_1^2(\Gamma,
\alpha_+\mu^{-n})$ with a suitable $n$. Therefore, 
linear
combinations of functions from 
$\{b^n(t)K^{\alpha_+\mu^{-n}}_{s_+b^{2n}}(t)\}$ are 
dense
in $L^2_{s_+}$. Since this system of functions
is orthonormal, it forms a basis in $L^2_{s_+}$.

Since $bz\in H^\infty(\Gamma,\mu)$, we have
$$
z:b^n H^2_{s_+b^{2n}}(\alpha_+\mu^{-n})\to
b^{n-1} H^2_{s_+b^{2n-2}}(\alpha_+\mu^{-n+1}).
$$
For this reason, in the basis
$\{b^n(t)K^{\alpha_+\mu^{-n}}_{s_+b^{2n}}(t)\}_{n\in\Bbb Z}$,
the matrix of the operator multiplication by $z(t)$ has
only one non--zero entry over diagonal in each column.
But the operator is self--adjoint, therefore, the matrix
is a three--diagonal Jacobi matrix.\qed
\enddemo

\proclaim{Lemma 1.5} Let
$e^+(n,t)=b^n(t)K^{\alpha_+\mu^{-n}}_{s_+b^{2n}}(t)$,
$n\in \Bbb Z$. Define
$$
s(t)e^-(n,t)=\bar t e^+(-n-1,\bar t)+
s_+(t) e^+(-n-1,t).
$$
Then $\{e^-(n,t)\}$ is an orthonormal basis in
$L^2_{s_-}$,
$$
s(t)e^-(n,t)\in A_1^2(\Gamma,\alpha^{-1}_+),\quad
n\in\Bbb Z_+,\tag 1.9
$$
and
$$
e^-(0,0) (b e^+)(-1,0)=\frac{b'(0)}{s(0)}.\tag 1.10
$$
\endproclaim

\demo{Proof} Lemma 1.3 and Lemma 1.4 imply immediately 
that $\{e^-(n,t)\}$ is an orthonormal basis in
$L^2_{s_-}$. Moreover, $s(t)e^-(n,t)\in L^2_{dm|\Bbb E}$.
To prove (1.9) consider a scalar product 
($f\in A_1^2(\Gamma,\alpha_+)$)
$$\align
\langle \bar t f(\bar t), s(t)e^-(n,t)\rangle=&
\frac 1 2\left\langle
\bmatrix f(t)\\ \bar t f(\bar t)\endbmatrix,
\bmatrix 1 & \overline{s_+(t)}\\
       s_+(t)& 1\endbmatrix
\bmatrix e^-(n,t)\\ \bar t e^-(n,\bar t)\endbmatrix
\right\rangle\\
 =&\frac 1 2\left\langle
\bmatrix f(t)\\ \bar t f(\bar
t)\endbmatrix,
\bmatrix 1 & \overline{s_+(t)}\\
       s_+(t)& 1\endbmatrix
\bmatrix (b^{-n-1}K^{\alpha_+\mu^{n+1}}_{s_+b^{-2n-2}})
(t)\\
\bar t
(b^{-n-1} K^{\alpha_+\mu^{n+1}}_{s_+b^{-2n-2}})(\bar
t)\endbmatrix
\right\rangle\\
=&\langle b^{n+1} f, K^{\alpha_+\mu^{n+1}}_{s_+b^{-2n-2}}
\rangle_{s_+b^{-2n-2}}=0,\quad\forall n\ge 0.
\endalign
$$

To prove (1.10), we write
$$
s(0)e^-(0,0)=\langle s(t)e^-(0,t), k^{\alpha_+^{-1}}(t)
\rangle.\tag 1.11
$$
Due to the Direct Cauchy Theorem, the reproducing kernel
$k^\alpha$ possesses the following property:
$$
\bar t k^{\alpha^{-1}_+}(\bar t)=\frac {b'(0)}{
k^{\alpha_+\mu}(0)}\frac{k^{\alpha_+\mu}(t)}{b(t)}.
\tag 1.12
$$
Substituting (1.12) in (1.11), we obtain
$$
\align
s(0)e^-(0,0)=&
\frac {b'(0)}{2 k^{\alpha_+\mu}(0)}
\left\langle
\bmatrix 1 & \overline{s_+(t)}\\
       s_+(t)& 1\endbmatrix
\bmatrix e^-(-1,t)\\ \bar t e^-(-1,\bar t)\endbmatrix,
\bmatrix (b^{-1}k^{\alpha_+\mu})(t)\\ \bar t 
(b^{-1}k^{\alpha_+\mu})(\bar t)\endbmatrix
\right\rangle\\
=&
\frac {b'(0)}{ k^{\alpha_+\mu}(0)}
\langle
 K^{\alpha_+\mu}_{s_+b^{-2}}(t),k^{\alpha_+\mu}(t)
 \rangle_{s_+b^{-2}}.
\endalign
$$
Using (1.3), we have
$$
\align
s(0)e^-(0,0)=&
\frac {b'(0)}{ k^{\alpha_+\mu}(0)
K^{\alpha_+\mu}_{s_+b^{-2}}(0)}
\lim_{\epsilon\to 0}
\langle
(\epsilon+I+\Cal H_{s_+b^{-2}})^{-1}k^{\alpha_+\mu},
k^{\alpha_+\mu}\rangle_{s_+b^{-2}}\\
=&
\frac {b'(0)}{ k^{\alpha_+\mu}(0)
K^{\alpha_+\mu}_{s_+b^{-2}}(0)}
\lim_{\epsilon\to 0}
\langle
(I+\Cal H_{s_+b^{-2}})
(\epsilon+I+\Cal H_{s_+b^{-2}})^{-1}k^{\alpha_+\mu},
k^{\alpha_+\mu}\rangle\\
=&
\frac {b'(0)}{ k^{\alpha_+\mu}(0)
K^{\alpha_+\mu}_{s_+b^{-2}}(0)}\{ k^{\alpha_+\mu}(0)-
\lim_{\epsilon\to 0}
\epsilon\langle
(\epsilon+I+\Cal H_{s_+b^{-2}})^{-1}k^{\alpha_+\mu},
k^{\alpha_+\mu}\rangle\}.
\endalign
$$
Since the limit (1.7) exists, finally, we get
$$
s(0)e^-(0,0)=
\frac {b'(0)}{ 
K^{\alpha_+\mu}_{s_+b^{-2}}(0)}=
\frac {b'(0)}{ (be^+)(-1,0)}.
$$
The lemma is proved.\qed
\enddemo

\proclaim{Lemma 1.6} Let $||s_+||<1$. Then
$$
K_{s_\pm}^{\alpha_\pm}(0)
K_{s_\mp b^{-2}}^{\alpha_\mp\mu}(0)=\frac
{b'(0)}{s(0)}.
$$
\endproclaim

\demo{Proof} Note, that operators 
$(I+\Cal H_{s_\pm b^n})$
are invertible. 

We use the notation of Lemma 1.5. As we know, 
$s(t)e^-(0,t)\in A_1^2(\Gamma,\alpha^{-1}_+)$. But, in
the case under consideration, $1/s\in
H^\infty(\Gamma,\alpha_+\alpha_-)$. Hence, the function
$e^-(0,t)$ itself belongs to $ A_1^2(\Gamma,\alpha_-)$.
Therefore, we can project each term onto 
$A_1^2(\Gamma,\alpha_-)$ in the relation
$$
\bar t (se^+)(-1,\bar t)=e^-(0,t)+\bar t (s_-e^-)(0,\bar t).
$$
On the right hand side we get
$$
P_{ A_1^2(\Gamma,\alpha_-)}
\{e^-(0,t)+\bar t (s_-e^-)(0,\bar t)\}
=(I+\Cal H_{s_-})e^-(0,t).
$$
To evaluate the left hand side, using (1.10), we write
$$\align
s(t)e^+(-1,t)
=&s(0)(be^+)(-1,0)\frac{k^{\alpha_-^{-1}\mu}(t)} {b(t)
k^{\alpha_-^{-1}\mu}(0)}+g(t)\\
=&\frac{b'(0)}{e^-(0,0)}\frac{k^{\alpha_-^{-1}\mu}(t)} 
{b(t)
k^{\alpha_-^{-1}\mu}(0)}+g(t),\quad g\in
 A_1^2(\Gamma,\alpha_-^{-1})
\endalign
$$
Using (1.12), we get
$$
P_{ A_1^2(\Gamma,\alpha_-)}
\{\bar t (s e^+)(-1,\bar t)\}=
\frac{k^{\alpha_-}(t)}{e^-(0,0)}=
(I+\Cal H_{s_-})e^-(0,t).
$$
Thus,
$$
e^-(0,t){e^-(0,0)}=(I+\Cal H_{s_-})^{-1}{k^{\alpha_-}}.
$$
In particular, $e^-(0,0)=K_{s_-}^{\alpha_-}(0)$,
and (1.10) becomes the statement of the lemma.\qed

\enddemo

\proclaim{Lemma 1.7}
Assume that for some Jacobi matrix $J$
there exists 
  a  pair of unitary transforms
$$
\Cal F^\pm:l^2(\Bbb Z)\to L^2_{s_\pm},
 \quad
(\Cal F^{\pm} Jf)(t)=z(t)(\Cal F^\pm f)(t),
$$
determining each other by the relations
$$
s(t)(\Cal F^\pm f)(t)=
\bar t (\Cal F^\mp f)(\bar t)+
s_\mp(t)(\Cal F^\mp f)(t),
$$
such that
$$
s\Cal F^\pm(l^2(\Bbb Z_\pm))\subset
A^2_1(\Gamma,\alpha^{-1}_\mp).
\tag 1.13
$$
As before, we put
$$
e^+(n,t)=(\Cal F^+ e_n)(t),\quad
e^-(n,t)=(\Cal F^- e_{-n-1})(t).\tag 1.14
$$
Then $e^\pm(n,t)$ has at the origin zero (poles)
of multiplicity $n,\ n>0$ ($-n,\ n<0$). Furthermore,
$\Cal
F^\pm(l^2(\Bbb Z_\pm))\supset H^2_{s_\pm}(\alpha_\pm)$, 
and,
hence,
 $$
e^\pm(0,0)\ge K^{\alpha_\pm}_{s_\pm}(0).\tag 1.15
$$
The equality in (1.15) takes place if and only if
$e^\pm(0,t)= K^{\alpha_\pm}_{s_\pm}(t)$.

\endproclaim

\demo{Proof} Let us show that the annihilator of 
the linear
space $A_1^2(\Gamma,\alpha^+)\subset L^2_{s_+}$ contains
$\Cal F^+\{l^2(\Bbb Z_-)\}$. For 
$f\in A_1^2(\Gamma,\alpha^+)$ and $e^+(-n-1,t),\ n\ge 0$,
we have
$$\align
\langle f( t), e^+(-n-1,t)\rangle_{s_+}
=&
\frac 1 2\left\langle
\bmatrix f(t)\\ \bar t f(\bar t)\endbmatrix,
\bmatrix 1 & \overline{s_+(t)}\\
       s_+(t)& 1\endbmatrix
\bmatrix e^+(-n-1,t)\\ \bar t e^+(-n-1,\bar t)
\endbmatrix
\right\rangle\\
=&
\langle  f( t), e^+(-n-1,t)+
\overline{ t s_+(t)}e^+(-n-1,\bar t)
\rangle\\ 
=&
\langle f( t),\bar t (s e^-)(n,\bar t)\rangle.
\endalign
$$
 By (1.13) and (DCT), the last scalar product equals 
zero.
Therefore,
$$
H^2_{s_+}(\alpha_+)=\text{clos}_{
L^2_{s_+}}A_1^2(\Gamma,\alpha^+)\subset
\{\Cal F^+(l^2(\Bbb Z_-))\}^{\bot}=
\Cal F^+(l^2(\Bbb Z_+)).
$$

Now, from the three-term recurrent relation
$$
z(t) s(t) e^+(n,t)=p_n s(t) e^+(n-1,t)+
q_n s(t) e^+(n,t)+
p_{n+1} s(t) e^+(n+1,t),
\tag 1.16
$$
and (1.13) it follows that $e^+(n,t)$, $n>0$, has
in the origin zero, at least of multiplicity $n$.

Since $ K^{\alpha_+}_{s_+}(t)\in
\Cal F^+(l^2(\Bbb Z_+))$, it possesses the decomposition
$$
 K^{\alpha_+}_{s_+}(t)=\sum_{n\in\Bbb Z_+}a_n e^+(n,t).
$$
Since $e^+(n,0)=0$, $n>0$, 
$$
a_0=\frac{K^{\alpha_+}_{s_+}(0)}{e^+(0,0)}
$$
in this decomposition. But,
$$
|a_0|^2\le\sum |a_n|^2= 
|| K^{\alpha_+}_{s_+}(t)||_{s_+}^2=1.
$$
Thus, (1.15) and the lemma are proved. \qed
\enddemo

\proclaim{Lemma 1.8 [8]}
Let 
$f\in L^\infty(\alpha^{-2})$.
Then 
$$
P_{A_1^2(\Gamma,\alpha)}\left\{\bar t (f b^n
K^{\alpha\mu^{-n}})(\bar t)\right\}\to 0,\quad n\to
+\infty, 
$$
where $P_{A_1^2(\Gamma,\alpha)} $ is the orthogonal
projection from $L^2_{dm\vert\Bbb E} $ onto
${A_1^2(\Gamma,\alpha})$.
 \endproclaim

\demo{Proof} Let us denote by $\Delta^\beta(t)$ an 
extremal
function of the  problem
$$
\Delta^\beta(0)=\sup\{\phi(0):\ 
\phi\in H^\infty(\Gamma,\beta),\ ||\phi||\le 1\}.
$$
Using properties 1), 2) of a group of Widom type
with (DCT), Theorem [5], and compactness of $\Gamma^*$,
for any $\epsilon>0$, 
we can find a finite covering of $\Gamma^*$
$$
\Gamma^*=\bigcup_{j=1}^{l(\epsilon)}\{\beta:\ 
\text{dist}(\beta,\beta_j)\le
\eta(\epsilon)\}
$$
such that
$$
2\left|1-\Delta^{\beta_j^{-1}\beta}(0)\frac
{K^{\beta_j}(0)}{K^{\beta}(0)}\right|\le\epsilon^2,
\quad
\text{dist}(\beta,\beta_j)\le\eta(\epsilon).
$$
 It means that
$$
\Vert(\Delta^{\beta_j^{-1}\beta}K^{\beta_j}) - 
K^\beta\Vert^2
\le 1+1-2\Delta^{\beta_j^{-1}\beta}(0)\frac
{K^{\beta_j}(0)}{K^{\beta}(0)}\le
\epsilon^2,
\quad
\text{dist}(\beta,\beta_j)\le\eta(\epsilon).
$$

For fixed $\beta$  one can find $n_0$ such that
$$
\Vert P_{b^n A_1^2(\Gamma,\alpha^2\beta^{-1}\mu^{-n})}
\bar t
(f K^\beta)(\bar t)\Vert\le
\epsilon,\ 
\forall n>n_0,
$$
 Therefore, there exists $n_0$ such that
$$
\Vert P_{b^n A_1^2(\Gamma,\alpha^2\beta_j^{-1}\mu^{-n})}
\bar t
(f K^{\beta_j})(\bar t)\Vert
 \le
\epsilon,\ 
\forall n>n_0,\ 1\le j\le l(\epsilon).
$$

Now, let $n>n_0=n_0(\epsilon)$ and let $\beta_j:$ 
$\text{dist}(\beta_j,\alpha\mu^{-n})\le\eta(\epsilon)$. 
For $h \in A^2_1(\Gamma,\alpha)$, we write
$$
\langle \bar t (f b^n K^{\alpha\mu^{-n}})(\bar t),  
h\rangle=
\langle \bar t (b^n f
[K^{\alpha\mu^{-n}}-\Delta^{\alpha\mu^{-n}\beta_j^{-1}}
K^{\beta_j}])(\bar t),  h\rangle +
\langle \bar t (b^n \Delta^{\alpha\mu^{-n}\beta_j^{-1}} 
f
K^{\beta_j})(\bar t),  h\rangle. 
$$
Then
 $$
\vert 
\langle \bar t (b^n f
[K^{\alpha\mu^{-n}}-\Delta^{\alpha\mu^{-n}\beta_j^{-1}}
K^{\beta_j}])(\bar t),  h\rangle
\vert
\le \Vert f\Vert\, \Vert h\Vert\,
\Vert K^{\alpha\mu^{-n}}-
\Delta^{\alpha\mu^{-n}\beta_j^{-1}}
K^{\beta_j}\Vert\le \epsilon||f||\,\Vert h\Vert,
$$
and
$$\align
|\langle \bar t (b^n \Delta^{\alpha\mu^{-n}\beta_j^{-1}} f
K^{\beta_j})(\bar t),  h\rangle|=&
|\langle \bar t ( f
K^{\beta_j})(\bar t), b^n\overline{
\Delta^{\alpha\mu^{-n}\beta_j^{-1}}(\bar t)} h\rangle|\\
\le&
\Vert P_{b^n A_1^2(\Gamma,\alpha^2\beta_j^{-1}\mu^{-n})}
\bar t
(f K^{\beta_j})(\bar t)\Vert\,||h||\le\epsilon||h||.
\endalign
$$ 
Therefore,
$$
|\langle P_{A_1^2(\Gamma,\alpha)}\left\{\bar t (f b^n
K^{\alpha\mu^{-n}})(\bar t)\right\}, h\rangle|\le
\epsilon(1+||f||)||h||.
$$
Putting $h=P_{A_1^2(\Gamma,\alpha)}\left\{\bar t (f b^n
K^{\alpha\mu^{-n}})(\bar t)\right\}$, we get
$$
||P_{A_1^2(\Gamma,\alpha)}\left\{\bar t (f b^n
K^{\alpha\mu^{-n}})(\bar t)\right\}||
\le\epsilon(1+||f||).
$$
The lemma is proved.\qed

\enddemo

\proclaim{Proposition 1.1}
Assume that for some Jacobi matrix $J$
there exists 
  a  pair of unitary transforms
$$
\Cal F^\pm:l^2(\Bbb Z)\to L^2_{s_\pm},
 \quad
(\Cal F^{\pm} Jf)(t)=z(t)(\Cal F^\pm f)(t),
$$
determining each other by the relations
$$
s(t)(\Cal F^\pm f)(t)=
\bar t (\Cal F^\mp f)(\bar t)+
s_\mp(t)(\Cal F^\mp f)(t),
\tag 1.17
$$
such that (1.13) holds.
 Then the following relations are equivalent:
$$\align
&e^\pm(n,t)=b^n(t) K^{\alpha_\pm\mu^{-n}}+o(1)
\quad \text{in}\  L^2_{s_+};
\tag 1.18\\
 &\bar t p_n\{e^+(n,t)e^+(n-1,\bar t)-
e^+(n-1,t)e^+(n,\bar t)\}=z'(t);
\tag 1.19\\
 &
{s(0)}e^+(0,0) (b e^-)(-1,0)={b'(0)},
\tag 1.20
\endalign
$$
 where $\{e^\pm(n,t)\}$ is defined by (1.14).
\endproclaim

\demo{Proof} (1.18) $\Rightarrow$ (1.19). It follows
from two remarks. First, the form on the left in
(1.19) does not depend on $n$ ( it is the Wronskian
of recurrence relation (0.7)). Second, the identity
$$
\bar t\frac{K^\alpha(0)}{K^{\alpha\mu}(0)}
\{K^\alpha(t)(K^{\alpha\mu}/b)(\bar t)-
(K^{\alpha\mu}/b)(t)K^\alpha(\bar t)\}=z'(t)
$$
holds for any $\alpha\in\Gamma^*$.

\medskip

(1.19) $\Rightarrow$ (1.20). Let us introduce the
matrix
$$
\Phi (t)=\bmatrix
e^-(-1,t)&-e^-(0,t)\\
-e^+(0,t)&e^+(-1,t)
\endbmatrix. \tag 1.21
$$
Then (1.17) implies
$$
\bar t \Phi(\bar t)=-S(t)\Phi(t).
$$
In particular, with a help of (1.19), we get
$$\align
s(t)=&-\bar t\frac{
e^+(0,t)e^+(-1,\bar t)-
e^+(-1,t)e^+(0,\bar t)}{
e^-(-1,t)e^+(-1, t)-
e^-(0,t)e^+(0, t)}\\
 =&\frac{ - z'(t)}{p_0
\{e^-(-1,t)e^+(-1, t)-
e^-(0,t)e^+(0, t)\}}.\tag 1.22
\endalign
$$
Since $b(t)e^\pm(-1,t)$ are holomorphic functions
(in fact, of Smirnov class)
$$
s(0)=\frac{ b'(0)}{p_0
(be^-)(-1,0)(b e^+)(-1, 0)}.
$$
Now, we only have to mention that
$p_0(be^\pm)(-1,0)=e^\pm(0,0)$.

\medskip
(1.20) $\Rightarrow$ (1.18). This is non--trivial
part of the proposition. The main step is to prove
that
$$
\lim_{n\to+\infty}\frac{(b^{-n}e^+)(n,0)}
{K^{\alpha_+\mu^{-n}}(0)}=1.
\tag 1.23
$$
By Lemma (1.7) we have an estimate from below
$$\align
(b^{-n}e^+)(n,0)&\ge
K^{\alpha_+\mu^{-n}}_{s_+b^{2n}}(0)\ge
\{(\epsilon+I+\Cal H_{s_+b^{2n}})^{-1}
k^{\alpha_+\mu^{-n}}\}
^{1/2}(0)\\
&=
\frac 1{\sqrt{1+\epsilon}} 
K^{\alpha_+\mu^{-n}}_{\frac{s_+}{1+\epsilon}b^{2n}}(0).
\tag 1.24
\endalign
$$
To get an estimate from above we use (1.20).

Let us note that due to the recurrence relation, the 
form
$$
p_n
\{e^+(n-1,t)e^-(-n-1, t)-
e^+(n,t)e^-(-n, t)\}
$$
also does not depend on $n$.
Thus, a relation like (1.20) holds for all $n$:
$$\align
(b^{-n}e^+)(n,0)(b^{n+1}e^-)(-n-1, 0)=&
p_n(b^{-n+1}e^+)(n-1,0)(b^{n+1}e^-)(-n-1, 0)\\
=&e^+(0,0)(be^-)(-1, 0)=b'(0)/s(0).
\endalign
$$
Therefore,
$$\align
(b^{-n}e^+)(n,0)=&
\frac{b'(0)}{s(0)}\frac{1}{(b^{n+1}e^-)(-n-1, 0)}\\
\le&
\frac{b'(0)}{s(0)}\frac{1}
{K^{\alpha_-\mu^{n+1}}_{s_-b^{-2n-2}}(0)}\\
\le&
\frac{b'(0)}{s(0)}\frac{1}{
\{(\epsilon+I+\Cal H_{s_-b^{-2n-2}})^{-1}
k^{\alpha_-\mu^{n+1}}\}
^{1/2}(0)}\\
=&
\frac{b'(0)}{s(0)}\frac{\sqrt{1+\epsilon}}
{K^{\alpha_-\mu^{n+1}}_{s_{\epsilon,-}b^{-2n-2}}(0)},
\tag 1.25\endalign
$$
where $s_{\epsilon,-}:=s_-/(1+\epsilon)$.

With the function
$s_{\epsilon,-}$, let us associate the functions
$s_{\epsilon}$, $s_{\epsilon,+}$ and the character
 $\alpha_{\epsilon,+}$ (note, that $s_{\epsilon,+}$
is not $\frac 1{1+\epsilon}s_+$, but
$s_{\epsilon,+}=-\bar s_{\epsilon,-} (s_\epsilon /
\bar s_{\epsilon})$). It is important that 
$s_\epsilon(0)$
and $\alpha_{\epsilon,+}$ depend continuously on
$\epsilon$.

By Lemma 1.6
$$
\frac{b'(0)}
{K^{\alpha_-\mu^{n+1}}_{s_{\epsilon,-}b^{-2n-2}}(0)}
= s_\epsilon(0)
{K^{\alpha_{\epsilon,+}
\mu^{-n}}_{s_{\epsilon,+}b^{2n}}(0)}.
\tag 1.26
$$
Substituting (1.26) in (1.25), and combining the result
with (1.24), we obtain
$$
\frac 1{\sqrt{1+\epsilon}} 
K^{\alpha_+\mu^{-n}}_{\frac{s_+}{1+\epsilon}b^{2n}}(0)
\le
(b^{-n}e^+)(n,0)\le
\sqrt{1+\epsilon}\frac{s_\epsilon(0)}{s(0)}
{K^{\alpha_{\epsilon,+}
\mu^{-n}}_{s_{\epsilon,+}b^{2n}}(0)}.
\tag 1.27
$$

Lemma 1.8 implies that for any 
$f\in L^\infty(\Gamma,\alpha_+^{-2})$ with $||f||<1$
we have
$$
\lim_{n\to+\infty}\frac{
K^{\alpha_+\mu^{-n}}_{f b^{2n}}(0)}
{K^{\alpha_+\mu^{-n}}(0)}\to 1.
$$
Indeed, 
$$
\align
|k^{\alpha_+\mu^{-n}}_{f b^{2n}}(0)
-k^{\alpha_+\mu^{-n}}(0)|
=&|\langle
\Cal H_{f b^{2n}}k^{\alpha_+\mu^{-n}},
(I+\Cal H_{f b^{2n}})^{-1}k^{\alpha_+\mu^{-n}}
\rangle|\\
=&
|\langle\bar t(f b^n k^{\alpha_+\mu^{-n}})(\bar t),
b^n(I+\Cal H_{f b^{2n}})^{-1}k^{\alpha_+\mu^{-n}}
\rangle|\\
\le&
|| P_{A_1^2(\Gamma,\alpha_+)}\{
\bar t(f b^n k^{\alpha_+\mu^{-n}})(\bar t)\}||\,
||b^n(I+\Cal H_{f b^{2n}})^{-1}k^{\alpha_+\mu^{-n}}||\\
\le&
\frac 1{1-||f||}
|| P_{A_1^2(\Gamma,\alpha_+)}\{
\bar t(f b^n k^{\alpha_+\mu^{-n}})(\bar t)\}||\,
||k^{\alpha_+\mu^{-n}}||\to 0,
\endalign
$$
as $ n\to+\infty$.

Also, since $\alpha_{\epsilon,+}$ depends continuously
on $\epsilon$ and $K^{\alpha_+}(0)$ is continuous on 
a compact group $\Gamma^*$, for any $\delta>0$ we can 
choose
$\epsilon$ so small that
$$
\frac{K^{\alpha_{\epsilon,+}\mu^{-n}}(0)}
{K^{\alpha_+\mu^{-n}}(0)}\le 1+\delta,\quad\forall n.
$$

Thus, returning to (1.27), we obtain
$$
\frac 1{\sqrt{1+\epsilon}} 
 \le{\lim_{n\to\infty}\inf}
\frac{(b^{-n}e^+)(n,0)}
{K^{\alpha_+\mu^{-n}}(0)}
\le
\lim_{n\to\infty}\sup
\frac{(b^{-n}e^+)(n,0)}
{K^{\alpha_+\mu^{-n}}(0)}\le
\sqrt{1+\epsilon}\frac{s_\epsilon(0)}{s(0)}
 (1+\delta).
$$
Since $\epsilon$ and $\delta$ are arbitrary small,
(1.23) is proved.

Now we are in a position to prove (1.18). Consider 
the norm
of the difference
$$
||e^+(n,t)-b^n K^{\alpha_+\mu^{-n}}||^2_{s_+}=
1 +||b^n K^{\alpha_+\mu^{-n}}||^2_{s_+}
-2\langle e^+(n,t),b^n K^{\alpha_+\mu^{-n}}
\rangle_{s_+}.
$$
Since
$$
||b^n K^{\alpha_+\mu^{-n}}||^2_{s_+}=
1+\langle b^n K^{\alpha_+\mu^{-n}},\bar t
(s_+b^n K^{\alpha_+\mu^{-n}})(\bar t)
\rangle,
$$
using Lemma 1.8, we conclude that
$$
||b^n K^{\alpha_+\mu^{-n}}||^2_{s_+}\to 1,\ n\to+\infty.
$$
Let us evaluate the scalar product
$$
\align
\langle e^+(n,t),b^n K^{\alpha_+\mu^{-n}}
\rangle_{s_+}=&
\langle s e^-(-n-1,t),\bar t (b^n K^{\alpha_+\mu^{-n}})
(\bar t)
\rangle\\
=&
\langle s e^-(-n-1,t), b^{-n}b^{-1} 
K^{\alpha_+^{-1}\mu^{n+1}}
\rangle\\
=&
\frac{s(0)(b^{n+1}e^-)(-n-1,0)}{
K^{\alpha_+^{-1}\mu^{n+1}}(0)}\\
=&
\frac{K^{\alpha_+\mu^{-n}}(0)}
{(b^{-n}e^+)(n,0)}\to 1,\ n\to +\infty.
\endalign
$$
The proposition is proved.\qed
\enddemo

The following theorem shows that an arbitrary function
$s_+$, possessing (1.1), (1.2), is the reflection
coefficient of a Jacobi matrix of Szeg\"o class.

\proclaim{Theorem 1.1} Let a function 
$s_+\in L^\infty(\Gamma,\alpha_+^{-2}),
\ ||s_+||\le 1$, $\overline{s_+(\bar t)}=s_+(t)$,
be such that that
$\log(1-|s_+|^2)\in L^1$. Let an outer function
$s$, $s(0)>0$, and $s_-$ be associated to $s_+$ by the
relations
$$
|s|^2=1-|s_+|^2,\quad s_-=-\bar s_+ s/\bar s.
$$
Then the system of functions
$$
e^+(n,t)=b^n K^{\alpha_+\mu^{-n}}_{s_+b^{2n}}
$$
forms an orthonormal basis in $L^2_{s_+}$. The dual
system, defined by
$$
s(t) e^-(n,t)=\bar t e^+(-n-1,\bar t)+
s_+(t) e^+(-n-1,t),
$$
forms an orthonormal basis in $L^2_{s_-}$. The subspaces
of $L^2_{s_\pm}$, that formed by functions with 
vanishing
negative  Fourier coefficients with respect to these 
basses,
are spaces of holomorphic character--automorphic forms,
moreover,
$$
s f^\pm\in A_1^2(\Gamma,\alpha_\mp^{-1})\ \text{if}\ 
f^\pm\in\text{\rm clos}_{L^2_{s_\pm}}
\text{\rm span}\{e^\pm(n,t):\
n\ge 0\}.
$$
Further,
$$
e^\pm(n,t)=b^n K^{\alpha_\pm\mu^{-n}}+o(1)
\ \text{in}\ L^2_{s_\pm},
$$
and with respect to these basses 
the operator multiplication
by $z(t)$ is a Jacobi matrix J of Szeg\"o class.
\endproclaim

\demo{Proof} All statements, besides the last one, only
summaries results of Lemmas 1.4, 1.5 and 
Proposition 1.1.
To prove that $J$ is of Szeg\"o class we evaluate its
spectral density $\rho(x)$.

Using the definition of the resolvent matrix--function,
we get
$$
R(z)=
\bmatrix
\left\langle
(z(t)-z)^{-1} {e^+(-1, t)},e^+(-1, t)\right
\rangle_{s_+}&
\left\langle
(z(t)-z)^{-1} {e^+(0, t)},e^+(-1, t)
\right\rangle_{s_+}\\
\left\langle
(z(t)-z)^{-1} {e^+(-1, t)},e^+(0, t)\right\rangle_{s_+}&
\left\langle
(z(t)-z)^{-1} {e^+(0, t)},e^+(0, t)\right\rangle_{s_+}
\endbmatrix.
$$
Note, that if $f^\pm\in L^2_{s_\pm}$ are related
by $s(t)f^-(t)=\bar t f^+(\bar t) +s_+(t)f^+(t)$ then
$$
\left(\frac {f^+(t)}{z(t)-z}\right)^-=
\frac {f^-(t)}{z(t)-z}.
$$
Therefore, using Lemma 1.3, we have
 $$
R(z)= \frac 1 2\int_{\Bbb E}
\bmatrix
{e^+(-1, t)}&{e^+(0, t)}\\
{e^-(0, t)}&{e^-(-1, t)}
\endbmatrix^*
\bmatrix
{e^+(-1, t)}&{e^+(0, t)}\\
{e^-(0, t)}&{e^-(-1, t)}
\endbmatrix
\frac{|s(t)|^2\,dm}{z(t)-z},
$$
 and, substituting $s(t)$ from (1.22), we obtain
$$\align
R(z)=&
\frac 1 2\int_{\Bbb E}
\frac{
\bmatrix
{e^+(-1, t)}&{e^+(0, t)}\\
{e^-(0, t)}&{e^-(-1, t)}
\endbmatrix^*
\bmatrix
{e^+(-1, t)}&{e^+(0, t)}\\
{e^-(0, t)}&{e^-(-1, t)}
\endbmatrix}
{p^2_0
|e^-(-1,t)e^+(-1, t)-e^-(0,t)e^+(0, t)|^2}
\frac{|z'(t)|^2\,dm}{z(t)-z}\\
=&
\frac 1 2\int_{\Bbb E}
\frac{\tilde\Phi^{-1 *}(t)\tilde\Phi^{-1}(t)}
{z(t)-z}\frac{|z'(t)|^2\,|dt|}{2\pi p_0^2},
\endalign
$$
where
$$
\tilde\Phi(t)=
\bmatrix
{e^-(-1, t)}&-{e^+(0, t)}\\
-{e^-(0, t)}&{e^+(-1, t)}
\endbmatrix.\tag 1.28
$$
Thus,
$$
2\pi p_0^2\rho(z(t))=
{\tilde\Phi^{-1 *}(t)\tilde\Phi^{-1}(t)}|z'(t)|,
\tag 1.29
$$
and
$$
\det\{2\pi p_0\rho(z(t))\}=
\frac{|z'(t)|^2}{p_0^2|\det\tilde\Phi(t)|^2}=
|s(t)|^2.
$$
The theorem is proved.\qed
\enddemo

Let us note, by the way, that
$\Phi(t)$ (see (1.21)) and $\tilde\Phi(t)$ are related
by $\tilde\Phi(\bar t)=\Phi^*(t)$ and, besides (1.28),
$$
2\pi p_0^2\rho(z(t))=
\Phi^{-1}(t)\Phi^{-1 *}(t)|z'(t)|.\tag 1.30
$$

%%%%%%%%%%%%%%%%%
%%%%%%%%%%%%%%%%%
\head From spectral data to scattering data
\endhead

We start this section with the remark that
the spectral measure $d\sigma$ determines
a Jacobi matrix uniquely, but it is not
an arbitrary $2\times 2$ matrix--measure,
or, say, a real-valued (all entries are real) 
$2\times 2$ matrix--measure.

Indeed, one can represent $J$ as a two
dimensional perturbation of an orthogonal
sum of a pair of one--sided Jacobi matrices,
i.e.:
$$
J=\bmatrix
J_-&0\\
0& J_+
\endbmatrix+ p_0
\langle\ ,e_{-1}\rangle e_0
+ p_0
\langle\ ,e_{0}\rangle e_{-1},
$$
where $J_\pm=P_{l^2(\Bbb Z_\pm)} J|
{l^2(\Bbb Z_\pm)}$. This formula implies that
$$
R(z)=\bmatrix
r_-^{-1}(z)&p_0\\
p_0&r_+^{-1}(z)
\endbmatrix^{-1},
\tag 2.1
$$
where
$$\align
r_-(z)=&r(z,J_-)=
\langle(J_- -z)^{-1}e_{-1},e_{-1}\rangle
=\int\frac{d\sigma_-(x)}{x-z},\\
r_+(z)=&r(z,J_+)=
\langle(J_+ -z)^{-1}e_{0},e_{0}\rangle
=\int\frac{d\sigma_+(x)}{x-z}.
\endalign
$$
Thus, the real--valued matrix--measure
$d\sigma$ is determined by two
scalar measures $d\sigma_\pm$
(with the normalization $\int\,d\sigma_\pm
=1$) and a constant $p_0$.

In what follows $\hat f(x)\in L^2_{d\sigma}$
denotes the image of $f\in l^2(\Bbb Z)$ in
the spectral representation. Recall that
$$
\hat e_{-1}=\bmatrix 1\\0\endbmatrix,\quad
\hat e_{0}=\bmatrix 0\\1\endbmatrix
$$
and
$$
(\widehat {Jf})(x)=x\hat f(x).
$$

Let $\{P_n^\pm(z)\}$ be the orthonormal
 polynomials with respect to the (scalar)
measure
$d\sigma_\pm$ and 
$$
Q^\pm_n(z):=\int\frac
{P_n^\pm(x)-P_n^\pm(z)}{x-z}d\sigma_\pm(x)
$$
(so--called polynomials of the second kind).
In these terms 
$$\aligned
\hat e_n(x)=&\bmatrix-p_0 Q^+_n(x)\\
P_n^+(x)\endbmatrix,
\quad n\ge 0,\\
\hat e_{-n-1}(x)=&\bmatrix
P_n^-(x)\\ -p_0 Q^-_n(x)\endbmatrix,
\quad n\ge 0.
\endaligned\tag 2.2
$$
 
\medskip
Now, we prove Theorem 0.5.

\demo{Proof of Theorem 0.5, the uniqueness
part} The function 
$e^\pm(0,\zeta)/e^\pm(-1,\zeta)$ is
$\Gamma$--automor\-phic, thus it defines a
meromorphic function in 
$\bar\Bbb C\setminus E$,
$$
\tilde r_\pm(z(\zeta)):=-\frac
{e^\pm(0,\zeta)}{p_0 e^\pm(-1,\zeta)}.
$$
The recurrence relations implies that
$\tilde r_\pm(z)$ possesses the same
decomposition into a continued fraction
as $ r_\pm(z)$. Therefore,
$$
 r_\pm(z(\zeta))=-\frac
{e^\pm(0,\zeta)}{p_0 e^\pm(-1,\zeta)}.
\tag 2.3
$$
By Proposition 1.1 the asymptotic (1.18)
implies the identity (1.19). Using this
identity, we get ($t\in\Bbb T$)
$$
r_\pm(z(t))-\overline{r_\pm(z(t))}
=-p_0\frac{e^\pm(0,t)
\overline{e^\pm(-1,t)}-
e^\pm(-1,t)\overline{e^\pm(0,t)}}
{|p_0 e^\pm(-1,t)|^2}=
\frac{-tz'(t)}{|p_0 e^\pm(-1,t)|^2}
$$
This means that an outer part of the function
$ e^\pm(-1,\zeta)$ is determined uniquely. But
then (2.3) means that an outer part of
$e^\pm(0,\zeta)$ is determined uniquely,
and
since $b(\zeta)e^\pm(-1,\zeta)$
and $e^\pm(0,\zeta)$ are of Smirnov class,
these functions are determined up to a common
inner factor $\Delta_{\pm}(\zeta)$, i.e.,
$$
e^\pm(0,\zeta)=
\Delta_{\pm}(\zeta)\tilde e^\pm(0,\zeta)
\ \text{and}\ 
e^\pm(-1,\zeta)=
\Delta_{\pm}(\zeta)\tilde e^\pm(-1,\zeta),
\tag 2.4
$$
where the inner parts of
$\tilde e^\pm(0,\zeta)$,
$\tilde e^\pm(-1,\zeta)$ are relatively prime.

To show that $\Delta_{\pm}(\zeta)=1$
we use (0.23), (0.24). Since
$$\aligned
s(t)e^\mp(0,t)=&\bar t
e^\pm(-1,\bar t)+s_\pm(t)
e^\pm(-1,t),\\
s(t)e^\mp(-1,t)=&\bar t
e^\pm(0,\bar t)+s_\pm(t)
e^\pm(0,t),\endaligned
\tag 2.5
$$
we have
$$
s(t)\{e^\mp(-1,t)e^\pm(-1,t)-
e^\mp(0,t)e^\pm(0,t)\}=
\bar t\{e^\pm(0,\bar t)e^\pm(-1,t)-
e^\pm(-1,\bar t)e^\pm(0,t)\}.
$$
Substituting (2.4) and using the symmetry
$$
\overline{\tilde e^\pm(0,\bar t)}=
\tilde e^\pm(0,t),\quad
\overline{\tilde e^\pm(-1,\bar t)}=
\tilde e^\pm(-1,t),
$$
we obtain
$$\aligned
&s(t)b^2(t)\{e^\mp(-1,t)\tilde e^\pm(-1,t)-
e^\mp(0,t)\tilde e^\pm(0,t)\}\\
&\ \ \ =
\bar t\Delta_\pm(\bar t) b^2(t)
\{\overline{\tilde e^\pm(0,t)}
\tilde e^\pm(-1,t)-
\overline{\tilde e^\pm(-1, t)}
\tilde e^\pm(0,t)\}\\
&\ \ \ =-b^2(t)z'(t)\{p_0
\overline{\Delta_\pm(\bar t)}\}^{-1}.
\endaligned
$$
Since the first expression here is a function
of Smirnov class and $b^2 z'$ is an outer
function we conclude that $\Delta_\pm( t)$
is a constant.

Since
$$
\bar t\Phi(\bar t)=-S(t) \Phi (t)\tag 2.6
$$
with $\Phi(t)$ defined by (1.21), $S(t)$
is also determined in a unique way. At last,
by the recurrence relations we get the same
conclusion with respect to  all functions
$\{e^\pm(n,\zeta)\}$, not only for $n=-1,0$.
\qed
\enddemo

\demo{ Proof of Theorem 0.5, the existence
part}
The key instrument is the following theorem
[9]: if $r(z)$ is a meromorphic function in
$\bar \Bbb C\setminus E$ such that 
$\text{Im}r(z)/\text{Im}z\ge 0$ and poles
of $r(z(\zeta))$ satisfies the Blaschke
condition, then $r(z(\zeta))$ is a function
of bounded characteristic in $\Bbb D$ without
a singular component in the multiplicative
representation.

Let us show that poles of $r_\pm(z(\zeta))$
satisfies the Blaschke condition. Diagonal
entries $R_{-1,-1}(z)$ and $R_{0,0}(z)$ of the
resolvent matrix--function $R(z)$ are
holomorphic in $\bar\Bbb C\setminus E$. By the
theorem mentioned above they are functions of
bounded characteristic. In force of (2.1),
$$\align
-1/R_{-1,-1}(z)&=-1/r_-(z)+p_0^2 r_+(z),\\
-1/R_{0,0}(z)&=-1/r_+(z)+p_0^2 r_-(z).
\endalign
$$
This means that poles of $r_\pm$ are subsets of
poles of $1/R_{-1,-1}$ and $1/R_{0,0}$. Thus
$r_\pm(z(\zeta))$ are functions of bounded
characteristic.

Now, let us use the Szeg\"o condition
$\log\det\text{Im}\ R(z(t))\in L^1$. Since
$$
\det\text{Im}\ R^{-1}(z(t))=
|\det R^{-1}(z(t))|^2\det\text{Im}\
R(z(t)),
$$
using again (2.1), we have
$$
\log\text{Im}\ r_-^{-1}(z(t))+
\log\text{Im}\ r_+^{-1}(z(t)) =
\log\det\text{Im}\ R^{-1}(z(t))\in L^1.
$$
Therefore, each of the functions
$\log\text{Im}\ r_\pm(z(t))$ belongs to
$L^1$. Thus we can represent
$r_\pm(z)$ (uniquely) in the form
$$
r_\pm(z(\zeta))=
-\frac{e^\pm(0,\zeta)}{p_0e^\pm(-1,\zeta)},
$$
where ${e^\pm(0,\zeta)}$
and $b(\zeta){e^\pm(-1,\zeta)}$ are
functions of Smirnov class with coprime inner
parts (in fact, they are  Blaschke products)
such that
$$
\bar t
p_0\{{e^\pm(0,t)}{e^\pm(-1,\bar t)}-
{e^\pm(-1,t)}{e^\pm(0,\bar t)}\}=z'(t),
\tag 2.7
$$
and ${e^\pm(0,0)}>0$,
${(be^\pm)(-1,0)}>0$. Note that
$$
p_0=\frac{e^\pm(0,0)}{(be^\pm)(-1,0)}.
$$

As soon as the functions
${e^\pm(0,\zeta)}$ and ${e^\pm(-1,\zeta)}$
have been constructed
we are able to  introduce $S(t)$ and $\Cal
F^\pm$ in their terms.

First, let us write down an expression for the
resolvent  matrix--function:
$$\aligned
R(z(\zeta))=&\bmatrix
-p_0\frac{e^-(-1,\zeta)}{e^-(0,\zeta)}
&p_0\\
p_0&
-p_0\frac{e^+(-1,\zeta)}{e^+(0,\zeta)}
\endbmatrix^{-1}\\
=&-(p_0\Phi)^{-1}\Psi=-
\tilde\Psi(p_0\tilde\Phi)^{-1},
\endaligned\tag 2.8
$$
where $\Phi$ and $\tilde\Phi$ are as in
(1.21) and (1.28) respectively, and
$$
\Psi(\zeta)=\tilde\Psi(\zeta)=
\bmatrix{e^-(0,\zeta)}&0\\
0&{e^+(0,\zeta)}
\endbmatrix.
$$
Therefore,
$$
p_0^2\{R(z(t))-R^*(z(t))\}=
-tz'(t)\Phi^{-1}(t)\Phi^{-1 *}(t)=
-tz'(t)\tilde\Phi^{-1 *}(t)
\tilde\Phi^{-1}(t),\tag 2.9
$$
since (see (2.7))
$$
p_0\{\Psi\Phi^*-\Phi\Psi^*\}=
p_0\{\tilde\Phi^*\tilde\Psi-
\tilde\Psi^*\tilde\Phi\}=tz'.
$$

From (2.9) and $\tilde \Phi^*(t)=\Phi(\bar t)$
we get immediately that the matrix--function
$S(t)$ defined by (2.6) is unitary--valued.
Let us show that its element $s(\zeta)$ is an
outer function. In fact, we have to show
that the function
$b^2(\zeta) \det \Phi(\zeta)
$ is an outer function (see (1.22)).
To this end let us use the representation
for the diagonal entries of $R(z)$ (see (2.8))
$$\align
R_{-1,-1}(z(\zeta))=&
-\frac{{e^+(-1,\zeta)}{e^-(0,\zeta)}}
{p_0\det\Phi(\zeta)},\\
R_{0,0}(z(\zeta))=&
-\frac{{e^-(-1,\zeta)}{e^+(0,\zeta)}}
{p_0\det\Phi(\zeta)}.
\endalign
$$

Let $\Delta$ be an inner part of
$b^2(\zeta) \det \Phi(\zeta)$. Since
$R_{0,0}(z(\zeta))$ is of Smirnov class,
$\Delta$ is a divisor of
${{e^-(-1,\zeta)}{e^+(0,\zeta)}}$.
If $\Delta$ is not trivial, then it has
a non--trivial divisor $\Delta_1$ that is a
divisor of
one of these functions, say,
${e^-(-1,\zeta)}$. Since
${e^-(-1,\zeta)}$ and ${e^-(0,\zeta)}$ are
coprime (and $\Delta_1$ is a divisor of
$b^2(\zeta) \det \Phi(\zeta)$), the
$\Delta_1$ is a divisor of $ {e^+(0,\zeta)}$,
and, therefore, it is not a divisor of
${e^+(-1,\zeta)}$. Thus, $\Delta_1$
is not a divisor of the product
${e^+(-1,\zeta)}{e^-(0,\zeta)}$.
But this means that $R_{-1,-1}(z(\zeta))$
is not of Smirnov class. We arrive to a
contradiction, hence $\Delta$ is a constant.

\medskip
We define $\Cal F^\pm$ by the formulas
$$\aligned
(\Cal F^+ f)(t)=&
\bmatrix e^+(-1,t)& e^+(0,t)\endbmatrix
\hat f(z(t)),\\
(\Cal F^- f)(t)=&
\bmatrix e^-(0,t)& e^-(-1,t)\endbmatrix
\hat f(z(t)).
\endaligned\tag 2.10
$$
Evidently, $(\Cal F^\pm Jf)(t)=
z(t)(\Cal F^\pm f )(t)$ and by (2.6),
(0.23) are fulfilled.  Using the formula for the
spectral density $\rho(x)=\frac 1\pi
\text{Im}\ R(x)$ and (2.9), we have
$$
\int_E\hat f^*(x)\,\rho(x)dx\,\hat f(x)=
\frac 1 2 \int_{\Bbb E}
\hat f^*(z(t))(p_0\tilde\Phi)^{* -1}
(p_0\tilde\Phi)^{ -1}\hat f(z(t))\,|z'(t)|^2
dm(t).
$$
Since
$$
\tilde\Phi^{-1}(t)=\frac 1{\det \tilde\Phi(t)}
\bmatrix e^+(-1,t)& e^+(0,t)\\
e^-(0,t)& e^-(-1,t)\endbmatrix,
$$
we obtain
$$
\align
||f||^2=||\hat f||^2_{L^2_{d\sigma}}=&
\frac 1 2 \{||s\Cal F^+ f||^2+
||s\Cal F^- f||^2\}\\
=&||\Cal F^+ f||^2_{s_+}=
||\Cal F^- f||^2_{s_-}.
\endalign
$$
Thus $\Cal F^+$ is an isometry, and since
this map is invertible,
$$
\bmatrix\hat f_{-1}(z(t))\\
\hat f_{0}(z(t))
\endbmatrix=-\frac{p_0}{z'(t)}
\bmatrix \bar t e^+(0,\bar t)& -e^+(0,t)\\
-\bar t e^+(-1,\bar t)& e^+(-1,t)\endbmatrix
\bmatrix g(t)\\ \bar t g(\bar t)
\endbmatrix,
$$
where $g=\Cal F^+ f$, it is a unitary map.

Further,  using (2.2), for $n\ge 0$
we have
$$
e^+(n,\zeta)=\bmatrix e^+(-1,\zeta)&
e^+(0,\zeta)\endbmatrix
\bmatrix -p_0 Q^+_n(z(\zeta))\\
P^+_n(z(\zeta))\endbmatrix.
$$
Due to the well known properties of orthogonal
polynomials these functions have no
singularity at the origin and hence
they are functions of Smirnov class.
This easily implies (0.24).

At last, our maps possess properties
(1.19) (or (1.20)), in force of 
Proposition 1.1, (0.25) holds.
The theorem is proved.\qed
\enddemo

\proclaim{Theorem 2.1} Let 
$s_+\in L^\infty(\Gamma,\alpha^{-2}_+)$,
$||s_+||\le 1$, $\overline{s_+(\bar t)}=s_+(t)$,
 satisfy
$\log(1-|s_+|^2)\in L^1$. Then the reflection
coefficient $s_+$ determines a Jacobi matrix
of Szeg\"o class in a unique way if and only
if
$$
s(0)K^{\alpha_\pm}_{s_\pm}(0)
K^{\alpha_\mp\mu}_{s_\mp b^{-2}}(0)=b'(0).
\tag 2.11
$$
\endproclaim

\demo{Proof} Assume on the contrary that
$$
s(0)K^{\alpha_+}_{s_+}(0)
K^{\alpha_-\mu}_{s_- b^{-2}}(0)\not=b'(0).
\tag 2.12
$$
We construct two
Jacobi matrices. First, we consider the basis
$$
e^+(n,t)=b^n(t)
K^{\alpha_+\mu^{-n}}_{s_+b^{2n}}(t),\tag 2.13
$$
and by $J$ we denote the operator
multiplication by $z(t)$ in $L^2_{s_+}$
with respect to this basis
(Lemma 1.4). Then, starting
with the basis $\{b^n(t)
K^{\alpha_-\mu^{-n}}_{s_-b^{2n}}(t)\}$
in $L^2_{s_-}$, we introduce the basis
$$
s(t)\tilde e^+(-n-1,t)=\bar t
(b^n
K^{\alpha_-\mu^{-n}}_{s_-b^{2n}})(\bar t)
+ s_-(t) (b^n
K^{\alpha_-\mu^{-n}}_{s_-b^{2n}})(t).
\tag 2.14
$$
By $\tilde J$
we denote the operator
multiplication by $z(t)$ in $L^2_{s_+}$
with respect to $\{\tilde e^+(n,t)\}$.
By Lemma 1.5,
$$
s(0)\tilde e^+(0,0)
K^{\alpha_-\mu}_{s_-b^{-2}}(0)=b'(0).
$$
Thus (see (2.12)),
$\tilde e^+(0,0)\not = e^+(0,0)$.
Due to the uniqueness part of Theorem 0.5,
$\tilde J\not= J$. The "only if" part is
proved.

\medskip

Now, let (2.11) holds, and let $J$ be a Jacobi
matrix of Szeg\"o class and
$\Cal F^\pm$ be its representations in 
$L^2_{s_\pm}$. By Lemma 1.7,
$$
K^{\alpha_\pm}_{s_\pm}(0)\le
e^\pm(0,0)=\frac{b'(0)}{s(0)}
\frac 1{(be^\mp)(-1,0)}\le
\frac{b'(0)}{s(0)}\frac 1
{K^{\alpha_\mp\mu}_{s_\mp b^{-2}}(0)}.
$$
Then (2.11) implies that, in fact,
$ 
e^\pm(0,0)=K^{\alpha_\pm}_{s_\pm}(0)$ and
$(be^\mp)(-1,0)=
K^{\alpha_\mp\mu}_{s_\mp b^{-2}}(0)$,
thus, due to a conclusion of Lemma 1.7,
$$
e^\pm(0,t)=K^{\alpha_\pm}_{s_\pm}(t),\quad
e^\mp(-1,t)=
b^{-1}(t)K^{\alpha_\mp\mu}_{s_\mp b^{-2}}(t).
$$
Recall that these functions determine
the functions $r_\pm(z)$ and the coefficient
$p_0$ (see (2.3)), and they, in their turn,
determine $J$. The theorem is proved.
\enddemo

\proclaim{Corollary 2.1} Let $J$ be a
Jacobi matrix of Szeg\"o class
with a homogeneous spectrum $E$. Let
$\rho(x)$ be the density of its spectral
measure and $S(t)$ be its scattering
matrix--function. If
$$
\int_E\rho^{-1}(x)\,dx<\infty,\tag 2.15
$$
then there is no other Jacobi matrix of
Szeg\"o class with the same scattering
matrix--function $S(t)$.  
\endproclaim

\demo{Proof} By virtue of (1.29),
(2.15) is equivalent to
$$
\int_{\Bbb E}\tilde \Phi(t)\tilde\Phi^*(t)
\,dm<\infty,
$$
that is $e^\pm(0,t)$ and $e^\pm(-1,t)$
belong to $L^2_{dm|\Bbb E}$. Then words
by words repetition of arguments in
the proof of Lemma 1.6 gives us
$$
\aligned
(I+\Cal H_{s_\pm})e^\pm(0,t)e^\pm(0,0)=&
k^{\alpha_\pm}(t),\\
(I+\Cal H_{s_\pm b^{-2}})
(b e^\pm)(-1,t)(b e^\pm)(-1,0)=&
k^{\alpha_\pm\mu}(t).
\endaligned
$$
Thus, $e^\pm(0,0)=K^{\alpha_\pm}_{s_\pm}(0)$
and 
$(b e^\pm)(-1,0)=
K^{\alpha_\pm\mu}_{s_\pm b^{-2}}(0)
$.
Since, generally, 
$$
s(0)e^\pm(0,0)(b e^\pm)(-1,0)=b'(0),
$$
(2.11) holds, the corollary is proved.
\enddemo

 To finish this section we give an 
example of a scattering matrix--function, which
does not determine a Jacobi matrix of Szeg\"o
class. Moreover, in this example, the
associated operators $(I+\Cal H_{s_\pm})$ are
invertible.

\demo{Example} Let 
$v_\pm\in H^\infty(\Gamma)$, $||v_\pm||<1$,
$\overline{v_\pm(\bar t)}=v(t)$, 
$v_\pm(0)=0$. 
Define outer 
functions $u_\pm$, $u_\pm(0)>0$, by
$$
|u_\pm|^2+|v_\pm|^2=1.
$$
Then, we put
$$
s^0_\pm=-\bar v_\pm u_\pm/\bar u_\pm.
$$
At last,
$$
S(t)=\bmatrix s_-&s\\
s&s_+\endbmatrix=
\bmatrix s_-^0&0\\
0&s_+^0\endbmatrix+
\bmatrix u_-&0\\
0&u_+\endbmatrix
\Cal E\left\{ I-
\bmatrix v_-&0\\
0&v_+\endbmatrix\Cal E
\right\}^{-1}
\bmatrix u_-&0\\
0&u_+\endbmatrix,
$$
where 
$$
\Cal E=\bmatrix\frac{1+\Delta} 2&
\frac{1-\Delta} 2\\ \frac{1-\Delta} 2&
\frac{1+\Delta} 2
\endbmatrix,
$$
and $\Delta$ is an inner function from 
$H^\infty(\Gamma)$, 
$\overline{\Delta(\bar t)}=\Delta(t)$.

In this case $\Cal H_{s_\pm}=
\Cal  H_{s^0_\pm}$, since their 
symbols differ by functions from
$ H^\infty(\Gamma,\alpha^{-2}_\pm)$, and
therefore $(I+\Cal H_{s_\pm})$ are invertible.
From the other hand, the 
coefficient $s$ is of the form
$$
s=\frac{u_+ u_- (1-\Delta)/2}
{1-(v_++v_-)(1+\Delta)/2+v_+ v_-\Delta},
$$
and because of the factor $(1-\Delta)/2$, 
$1/s$ does not belong to $H^\infty(\Gamma,
\alpha_+\alpha_-)$.

The simplest choice of parameters: 
$$
E=[-2,2],\  \  v_\pm(t)=a_\pm t,
\ 0<a_\pm<1,
$$
$\Delta(t)$ is a Blaschke product,
$\deg\Delta>1$, gives us an
example where 
$\tilde e^+(-1,t)$, defined by (2.14), does not
belong to
$L^2$, at the same time $e^+(-1,t)$, defined by
(2.13), belongs to $ L^2$.
\enddemo   

%%%%%%%%%%%%%%%%%%%%%%
%%%%%%%%%%%%%%%%%%%%%
%%%%%%%%%%%%%%%%%%%%%
\head The Hilbert transform and $A_2$ condition
\endhead
By $\frak H$ we denote the transform
$$
(\frak H g)(z)=\int_{E}\frac{g(x)}{z-x}dx,
\tag 3.1
$$
primarily defined on  integrable 2D 
vector--functions.

\proclaim{Lemma 3.1} Let $J$ be of Szeg\"o
class and $\Cal F^\pm$ give its representations
as the operator multiplication by $z(t)$ in
the model spaces $L^2_{s_\pm}$. Then
$$
\bmatrix
\Cal F^- f^-\\ \Cal F^+ f^+
\endbmatrix(\zeta)=
p_0\Phi(\zeta)\{\frak H(\rho\hat
f)\}(z(\zeta)),
\tag 3.2
$$
for any finite vector
$f=f^-\oplus f^+\in 
l^2(\Bbb Z)=l^2(\Bbb Z_-)
\oplus l^2(\Bbb Z_+)$.
\endproclaim

\demo{Proof} Let $\tilde\Cal P_n(z)$ denote
the $n$--th matrix orthonormal polynomial
with respect the spectral measure $d\sigma$.
Recall, that
$$
\tilde\Cal P_n(z)=\bmatrix
\hat e_{-n-1}(z)&\hat e_n(z)\endbmatrix=
\bmatrix
P^-_n(z)&-p_0 Q_n^+(z)\\
-p_0 Q_n^-(z)&P^-_n(z)\endbmatrix,
$$
and, analogically to the scalar case,
$$
\tilde\Cal Q_n(z):=
\int d\sigma(x)\frac{\tilde\Cal P_n(z)-
\tilde\Cal P_n(x)}{z-x}=
\bmatrix
 Q_n^-(z)& 0\\
0& Q_n^+(z)\endbmatrix.
\tag 3.3
$$
Based on (3.3), we have
$$\align
p_0\Phi(\zeta)
\int d\sigma(x)\frac{
\tilde\Cal P_n(x)}{z(\zeta)-x}
=&p_0\Phi(\zeta)\left\{
\int \frac{d\sigma(x)}
{z(\zeta)-x}\tilde\Cal P_n(z(\zeta)) -
\int d\sigma(x)\frac{\tilde\Cal P_n(z(\zeta))-
\tilde\Cal P_n(x)}{z(\zeta)-x}
\right\}\\
=&
-p_0\Phi(\zeta)R(z(\zeta))
\tilde\Cal P_n(z(\zeta)) 
-p_0\Phi(\zeta)
\tilde\Cal Q_n(z(\zeta)).
\endalign
$$
Using (2.8) and definition (2.10), we get
$$\align
p_0\Phi(\zeta)
\int \frac{\rho(x)
\bmatrix
\hat e_{-n-1}(x)&\hat e_n(x)\endbmatrix
}{z(\zeta)-x}dx
=&
\Psi(\zeta)
\tilde\Cal P_n(z(\zeta))
-p_0\Phi(\zeta)
\tilde\Cal Q_n(z(\zeta)) \\
 =&\bmatrix
e^-(n,\zeta)&0\\
0&e^+(n,\zeta)\endbmatrix\\
=&\bmatrix
(\Cal F^- e_{-n-1})(\zeta)&0\\
0&(\Cal F^+ e_n)(\zeta)\endbmatrix.
\endalign
$$
In fact, this finishes the proof.\qed

\enddemo

\proclaim{Theorem 3.1} Let $\rho(x)$ be the 
spectral density of a Jacobi matrix $J$ of
Szeg\"o class and $s_+(t)$ be the  reflection
coefficient. Then the following statements
are equivalent:

1. There exist $C<\infty$ such that
$$\align
\int_E(\frak H g)^*(x-i0)\rho^{-1}(x)
(\frak H g)(x-i0)\,dx+&
\int_E(\frak H g)^*(x+i0)\rho^{-1}(x)
(\frak H g)(x+i0)\,dx\\
\le &C
\int_E g^*(x)\rho^{-1}(x) g(x)\,dx.
\tag 3.4
\endalign
$$

2. $s_+$ determines $J$ and the operators
$(I+\Cal H_{s_\pm})$ are invertible.
\endproclaim

\demo{Proof} $1\Rightarrow 2$. Since
(see (1.30))
$$\align
||\Cal F^-f^-||^2+||\Cal F^+f^+||^2
=&\int_{\Bbb E}
\{\frak H(\rho\hat f)\}^*(z(t))
(p_0\Phi)^*(t)(p_0\Phi)(t)
\{\frak H(\rho\hat f)\}(z(t))\,dm\\
=&\int_{\Bbb E}
\{\frak H(\rho\hat f)\}^*(z(t))
\frac 1 {2\pi}\rho^{-1}(z(t))
\{\frak H(\rho\hat f)\}(z(t))
\frac{|z'(t)||dt|}{2\pi}\\
=&\left(\frac 1 {2\pi}\right)^2\int_{ E}
\{\frak H(\rho\hat f)\}^*(x-i0)
\rho^{-1}(x)
\{\frak H(\rho\hat f)\}(x-i0)
\,dx\\
+&\left(\frac 1 {2\pi}\right)^2\int_{ E}
\{\frak H(\rho\hat f)\}^*(x+i0)
\rho^{-1}(x)
\{\frak H(\rho\hat f)\}(x+i0)
\,dx\\
\le&\frac C{(2\pi)^2}\int_{E}
\hat f^*(x)\rho(x) \hat f(x)\,dx=
\frac C{(2\pi)^2}||f||^2,\tag 3.5
\endalign
$$
we get $\Cal F^\pm f^\pm\in
A_1^2(\Gamma,\alpha_\pm)$. Thus,
$\Cal F^\pm\{l^2(Z_\pm)\}= 
H^2_{s_\pm}(\alpha_\pm)$. By Lemma 1.7
and Theorem 2.1
we come to the conclusion that
$s_+$ determines $J$. Further, by (3.5)
 $$\align
||\Cal F^-f^-||^2+||\Cal F^+f^+||^2
\le&
\frac C{(2\pi)^2}\{||f^-||^2+||f^+||^2\}\\
=&
\frac C{(2\pi)^2}\{
||\Cal F^-f^-||^2_{s_-}+
||\Cal F^+f^+||^2_{s_+}\}.
\endalign
$$
Using again $\Cal F^\pm f^\pm\in
A_1^2(\Gamma,\alpha_\pm)$, we can represent
the last norms in the form
$$
||\Cal F^-f^-||^2+||\Cal F^+f^+||^2
\le
\frac C{(2\pi)^2}\{
\langle(I+\Cal H_{s_-})\Cal F^-f^-,
\Cal F^-f^-\rangle+
\langle(I+\Cal H_{s_+})\Cal F^+f^+,
\Cal F^+f^+\rangle
\}.\tag 3.6
$$
This proves the second statement in 2.

$2 \Rightarrow 1$. Recall that
$H^2_{s_\pm}(\alpha_\pm)=\text{clos}_
{L^2_{s_\pm}} A_1^2(\Gamma,\alpha_\pm)$,
but in the case under consideration,
the norm in $H^2_{s_\pm}(\alpha_\pm)$ is
equivalent to the norm in 
$ A_1^2(\Gamma,\alpha_\pm)$, i.e.:
$$
h\in H^2_{s_\pm}(\alpha_\pm)
\Rightarrow h\in A_1^2(\Gamma,\alpha_\pm).
$$
Further, since $s_+$ determines $J$, by 
Lemma 1.7, we have 
$\Cal F^\pm\{l^2(Z_\pm)\}= 
H^2_{s_\pm}(\alpha_\pm)$. So, starting with
(3.6) we obtain (3.4).\qed
\enddemo

\subhead Remark on ISP and RHP
\endsubhead
Let us define $\sqrt{-b^2 z'}(\zeta)$ as the
square root of an outer function such that
$\sqrt{-b^2 z'}(0)>0$. Put
$$\sqrt{-z'}(\zeta)=
\frac{\sqrt{-b^2 z'}(\zeta)}{b(\zeta)}.
\tag 4.1
$$
In this case, 
$\overline{\sqrt{-z'}(\bar\zeta)}=
\sqrt{-z'}(\zeta)$. Let
$\Bbb E_-=\{t\in\Bbb E: \text{Im} t<0\}$.
Then $-t z'(t)=i|z'(t)|$, $t\in \Bbb E_-$.
Thus,
$$
t\{\sqrt{-z'}(t)\}^2=
i\sqrt{-z'}(t)\overline{\sqrt{-z'}(t)}
=i\sqrt{-z'}(t)\sqrt{-z'}(\bar t),
\quad t\in \Bbb E_-,
$$
or
$$
\bar t\sqrt{-z'}(\bar t)=-i\sqrt{-z'}(t),
\quad t\in \Bbb E_-.
\tag 4.2
$$
Besides,
$$
\sqrt{-z'}|[\gamma]=\epsilon(\gamma)
\sqrt{-z'}(\zeta),
$$
where $\epsilon\in\Gamma^*$,
$\epsilon^2=1_{\Gamma^*}$. But, in fact, the
group $\Gamma$ is defined up to a choice
of a half--period $\tilde\epsilon\in\Gamma^*$.
So, we may assume that 
$$
\sqrt{-z'}|[\gamma]=
\sqrt{-z'}(\zeta). \tag 4.3
$$

\proclaim{Lemma 4.1} Let
$E=[b_0,a_0]\setminus\cup_{j\ge 1}(a_j,b_j)$
be a homogeneous set. Then
$$
\Pi(z(\zeta)):=\frac{\Phi(\zeta)}
{\sqrt{-z'}(\zeta)}
$$
is a holomorphic matrix function
in $\bar\Bbb C\setminus[b_0,a_0]$ satisfying
the following RHP
$$\align
\Pi(x-i0)=&\bmatrix \alpha_{-,j}&0\\
0&\alpha_{+,j}\endbmatrix
\Pi(x+i0),\quad x\in(a_j,b_j),\tag 4.4\\
 \Pi(x-i0)=& -i\Sigma(x)
\Pi(x+i0),\quad x\in E,\tag 4.5
\endalign
$$
where $\Sigma(z(t)):= S(t)$, $t\in\Bbb E_-$,
with the normalization at infinity:
$$
\Pi(z)=
\bmatrix 1+\dots&-\frac a z+\dots\\
-\frac b z+\dots& 1+\dots
\endbmatrix
\bmatrix\frac 1{\sqrt{b s(0)}}&0\\
0&\frac 1{\sqrt{a s(0)}}
\endbmatrix.\tag 4.6
$$
\endproclaim

\demo{Proof} (4.4) follows from (4.3) and
(4.5) follows from (4.2) and (2.6). To prove
(4.6), we represent $\Pi(z)$ in the form
$$
\Pi(z)=
\bmatrix 1+\dots&-\frac a z+\dots\\
-\frac b z+\dots& 1+\dots
\endbmatrix
\bmatrix c_1&0\\
0&  c_2
\endbmatrix.
$$
Then, we note that
$$
c_1 c_2=\det \Pi(\infty)=
\left.\frac{e^-(-1,\zeta)e^+(-1,\zeta)-
e^-(0,\zeta)e^+(0,\zeta)}{-z'(\zeta)}
\right|_{\zeta=0}=\frac 1{p_0 s(0)},\tag 4.7
$$
and
$$
\frac{a c_2}{c_1}=\frac{b c_1}{c_2}=
\left.\frac{z(\zeta)e^\pm(0,\zeta)}
{e^\pm(-1,\zeta)}\right|_{\zeta=0}= p_0.
\tag 4.8
$$
Solving together (4.7), (4.8),
we get (4.6).\qed

\enddemo

\head Matrix $A_2$ on homogeneous sets
\endhead

Let $E$ be a homogeneous set. Throughout this
section $P_+$ denotes orthoprojector from the
vector--valued $L^2(\Bbb C^n)$ onto $H^2(\Bbb
C^n)$ in the upper halfplane. We are interested
in the boundedness of the weighted transform
$$
W^{1/2}P_+ W^{-1/2}:\chi_E L^2(\Bbb C^n)
\to\chi_E L^2(\Bbb C^n),
\tag 5.1
$$
 where $W$ is a weight
on $E$ and $\chi_E$ is the characteristic
function of the set $E$.

Here is an analog
of the matrix $A_2$ condition
$$
\sup_{x\in E, 0<\delta<1}||
\langle W\rangle^{1/2}_{I_{(x,\delta)}}
\langle W^{-1}\rangle^{1/2}_{I_{(x,\delta)}}||
<\infty,\tag 5.2
$$
where $I_{(x,\delta)}:=(x-\delta, x+\delta)$
and
$$
\langle W\rangle_{I_{(x,\delta)}}
:=\frac{1}{|I_{(x,\delta)}|}
\int_{I_{(x,\delta)}\cap E} W(t)\,dt.
$$
This supremum will be called $Q_{2,E}(W)$.

\proclaim
{Theorem 5.1} The operator (5.1) is bounded if
and only if $Q_{2,E}(W)<\infty$.
\endproclaim

\demo{Proof of necessity} With an arbitrary
$z_0\in \Bbb C_+$ we associate a subspace
$K_{b_{z_0}}=H^2(\Bbb C^n)\ominus b_{z_0}
H^2(\Bbb C^n)$ of the Hardy space, 
$b_{z_0}(z)=\frac{z-z_0}{z-\overline{z_0}}$.
It is well known, that
$$
P_{K_{b_{z_0}}}=P_+-b_{z_0}P_+\overline
{b_{z_0}}
$$
and
$$
\langle P_{K_{b_{z_0}}} f, g\rangle_
{L^2(\Bbb C^n)}=
\langle (P_+ f)(z_0), (P_+ g)(z_0)\rangle_
{C^n}.
$$
Because of the first of these relations
we have
$$
|\langle W^{1/2}P_{K_{b_{z_0}}} W^{-1/2}f,
g\rangle|\le 2Q||\chi_E f||\,||\chi_E g||.
$$
Now, using the second one we get
$$
|\langle (P_+ W^{-1/2}f)(z_0), (P_+
W^{1/2}g)(z_0)\rangle|\le 2Q||\chi_E
f||\,||\chi_E g||.
\tag 5.3
$$
Let us substitute 
$$
f=W^{-1/2}\frac{\xi}{x-\overline{z_0}},
\quad
g=W^{1/2}\frac{\eta}{x-\overline{z_0}},
\quad \xi, \ \eta\in \Bbb C^n
$$
in (5.3). This give us
$$
|\langle \,\langle W^{-1}\rangle_{z_0}
\xi, \langle W\rangle_{z_0}
\eta
\rangle|\le 2Q||
\langle W^{-1}\rangle^{1/2}_{z_0}
\xi||\,||\langle W\rangle^{1/2}_{z_0}
\eta||,
$$
where $\langle W\rangle_{z_0}$ denotes
an average with the Poisson kernel,
$$
\langle W\rangle_{z_0}:=
\frac 1{\pi }\int W\frac{\text{Im}\  
z_0}{|x-z_0|^2}\,dx.
$$
Thus we proved an inequality with the Poisson's
averages
$$
\langle W\rangle_{z_0}\le 2 Q
\langle W^{-1}\rangle_{z_0}^{-1}\tag 5.4
$$
At last let us note that
$$
\frac{\text{Im}\  
z_0}{|x-z_0|^2}\ge \frac{c}{|I|}
\chi_{I}, \quad
I=I_{(\text{Re}\ z_0,\text{Im}\ z_0)},
$$
 with an absolute and positive constant 
$c$. Therefore (5.4) implies (5.2).
\enddemo

\proclaim{Lemma 5.1} If $I$ is a centered at
$E$ interval and $z_0$ is the center of the
square built on $I$, then
$$
W\in A_2(E)\Rightarrow\quad
\langle W\rangle_{z_0}\le C(E, Q_{2,E}(W))
\langle W\rangle_{I}.\tag 5.5
$$
\endproclaim

\demo{Proof}
First we note, that for $\lambda=2/\eta$
$$
|\lambda I\cap E|\ge \eta|\lambda I|
\ge 2|I|\ge 2|I\cap E|,
$$
and therefore $|(\lambda I\setminus I)\cap E|
\ge |I\cap E|$. Let us show that 
$$
W(\lambda I)\ge \left(1+\frac{\eta^2}{\lambda^2 Q^2}
\right)W(I)
\quad\text{for}\ 
W\in A_2(E). \tag 5.6
$$
Integrating the inequality
$$
\bmatrix W^{-1}& 1\\
         1 & W\endbmatrix\ge 0
$$
over $(\lambda I\setminus I)\cap E$ we get
$$
\bmatrix W^{-1}(\lambda I\setminus I)& 
|(\lambda I\setminus I)\cap E|\\
 |(\lambda I\setminus I)\cap E| & 
W(\lambda I\setminus I)\endbmatrix\ge 0.
$$
Therefore
$$
\bmatrix W^{-1}(\lambda I)& 
|(\lambda I\setminus I)\cap E|\\
 |(\lambda I\setminus I)\cap E| & 
W(\lambda I)-W(I)\endbmatrix\ge 0,
$$
or
$$
W(\lambda I)-W(I)\ge|(\lambda I\setminus I)\cap E|^2
\{W^{-1}(\lambda I)\}^{-1}
\ge|I\cap E|^2
\{W^{-1}(\lambda I)\}^{-1}.
$$
Using (5.2) we obtain
$$
W(\lambda I)-W(I)\ge\frac{|I\cap E|^2}
{Q^2 |\lambda I|^2}W(\lambda I)\ge
\frac{\eta^2}{Q^2\lambda^2}W(I).
$$

To prove (5.5), using
$$
\frac{\text{Im}\  
z_0}{|x-z_0|^2}\le \frac{c}{|I|}
\sum\frac{1}{ \lambda^{2k}}
\chi_{\lambda^k I},\quad \lambda^k I=
I_{(\text{Re}\ z_0,\lambda^k\text{Im}\ z_0)}
$$ 
we write the following chain of inequalities:
$$
\align
\langle W\rangle_{z_0}&\le\frac{c}{|I|}
\sum\frac{1}{ \lambda^{2k}} W(\lambda^k I)\\
&\le\frac{c}{|I|}
\sum\frac{Q^2|\lambda^k I|^2}{ \lambda^{2k}}
\{W^{-1}(\lambda^k I)\}^{-1}\\
&\le\frac{c
Q^2}{|I|}\sum\left(1+\frac{\eta^2}{\lambda^2 Q^2}
\right)^{-k} |I|^2
\{W^{-1}(I)\}^{-1}\\
&\le\frac{c
Q^2}{\eta^2}\sum\left(1+\frac{\eta^2}{\lambda^2 Q^2}
\right)^{-k} 
\langle W\rangle_{I}.
\endalign
$$
\enddemo

\demo{Proof of sufficiency} We want to prove
that (5.2) suffices for $W^{1/2} P_+ W^{-1/2}$ in
(5.1) to be bounded. Fix $f, g\in \chi_E L^2(\Bbb C^n)$. 
We need to show
$$
\int_{\Bbb C_+}\left|
\langle(P_+ W^{-1/2}f)'(z), (P_+ W^{1/2}g)'(z)
\rangle_{\Bbb C^n}\right|\text{Im}\ z\,d A(z)\le
C||f||\,||g||.
$$
In other words, introducing a Stolz cone
$\Gamma_t$ and
$$
S(t)=
\int_{\Gamma_t}
\left|
\langle(P_+ W^{-1/2}f)'(z), (P_+ W^{1/2}g)'(z)
\rangle_{\Bbb C^n}\right|\,d A(z)
$$
one needs to prove that
$$
\int S(t)\, dt\le C||f||\,||g||.
\tag 5.7
$$
We follow closely the lines of the proof in [10].
Let us consider a nonnegative function $h(t)$
 and
$$
S_{h(t)}(t)=
\int_{\Gamma_{t, h(t)}}
\left|
\langle(P_+ W^{-1/2}f)'(z), (P_+ W^{1/2}g)'(z)
\rangle_{\Bbb C^n}\right|\,d A(z),
$$
where
$$
\Gamma_{t, h(t)}=\Gamma_t\cap\{z:\ \text {Im}\ z\le
h(t)\}.
$$
Let us note that
$$
\int S(t)\,dt\le c \int S_{h(t)}\, dt
\tag 5.8
$$
if the function $h(t)$ has the following property:
$$
\forall I\subset\Bbb R\quad
|\{t\in I:\ h(t)\ge |I|\}|\ge a|I|.
\tag 5.9
$$
Let us choose $h$ to be maximal such that
$$
S_{h(t)}(t)\le B(M||f||^{p_*})^{1/p_*}(t)
(M||g||^{p_*})^{1/p_*}(t),\tag 5.10
$$
where $B$, $p_*\in (1,2)$ will be chosen a bit
later and $M$ denotes the maximal function
$$
(M f)(x)=
\sup_{\delta>0}\frac{1}{|I(x,\delta)|}
\int_{I(x,\delta)} |f(t)|\,dt.
$$
If this $h$ satisfies (5.10), then (5.9) and
(5.8) imply what we need.

To choose $B, p_*$ and to prove that $h$ satisfies
(5.9) we follow the algorithm below.
Let $I_0$ be an arbitrary interval on the real axis.
We will consider two cases: $2I_0\cap E\not=\emptyset$
and $2I_0\cap E=\emptyset$.

In the first case we fix an interval $I$ centered at $E$
such that $I_0\subset I$ and $|I|\le 3|I_0|$.
Let $f_1=f\cdot\chi_{2I}$, $g_1=g\cdot\chi_{2I}$
and $f_2=f-f_1$, $g_2=g-g_1$. Denote
$A_I=\langle W\rangle_I^{1/2}$. Consider
$$
\aligned
t\in I,\quad &
S^{A_I}(f_i)(t)=\left(\int_{\Gamma_{t,|I|}}
||(P_+ A_I W^{-1/2}f_i)'||^2\,d A(z)\right)^{1/2},
\quad  i=1,2;\\
t\in I,\quad &
S^{A_I^{-1}}(g_i)(t)=\left(
\int_{\Gamma_{t,|I|}}
||(P_+ A_I^{-1} W^{1/2}g_i)'||^2\,d A(z)\right)^{1/2},
\quad  i=1,2.
\endaligned
$$
We will fix later $\alpha=\alpha(Q_{2,E}(W), n)>1$.
Now,
$$
\align
\frac{1}{|I|}&\int_{I}
\left(S^{A_I^{-1}}(g_1)(t)\right)^\alpha\,dt
\le 
\frac{C(\alpha)}{|I|}\int_{2I}||A_I^{-1} W^{1/2}g
(t)||^\alpha
\,dt\\
&
\le
\frac{ C_1(\alpha, n)}{|I|}\int_{2I}\left(||g
(t)||\sum_{i=1}^n ||W^{1/2}(t) A_I^{-1} e_i||
\right)^\alpha
\,dt\\
&
\le
\frac{ C_2(\alpha, n)}{|I|}\left(
\int_{2I} ||g(t)||^{(2-\tilde\epsilon)\alpha}
\,dt\right)^{\frac{1}{2-\tilde\epsilon}}
\sum_{i=1}^n \left(\int_{2I}
||W^{1/2}(t) A_I^{-1}
e_i||^{(2+\epsilon)\alpha}\,dt\right)^
{\frac{1}{2+\epsilon}}.
\tag 5.11\endalign
$$
Here $(2+\epsilon)^{-1}+(2-\tilde\epsilon)^{-1}=1$.
Notice that for every vector $x\in \Bbb C^n$
the scalar function $t\to ||W(t)^{1/2}x||$ is uniformly
in scalar $A_2(E)$. In particular, there exists such
an $\epsilon_0>0$ that we have the inverse H\"older
inequality for all such functions uniformly:
$$
\forall I\text{\ centered at\ } x\in E\quad
\left(\frac{1}{|I|}\int_I ||W(t)^{1/2}x||^{2+\epsilon_0}
\,dt\right)^{\frac{1}{2+\epsilon_0}}\le C
\left(\frac{1}{|I|}\int_I ||W(t)^{1/2}x||^{2}
\,dt\right)^{\frac{1}{2}}. \tag 5.12
$$
Let us choose $\epsilon=\frac{\epsilon_0}2$
$(\tilde\epsilon=\frac{\epsilon_0}{2+\epsilon_0})$,
$\alpha=1+\frac{\epsilon_0}{2(2+\epsilon_0)}$,
then   we have
$$
(2+\epsilon)\alpha<2+\epsilon_0,\tag 5.13
$$
$$
(2-\tilde\epsilon)\alpha<2.\tag 5.14
$$
We use (5.13) and inverse H\"older inequality (5.12)
in (5.11) to rewrite it as
$$
\align
&\left(\frac{1}{|I|}\int_{I}
\left(S^{A_I^{-1}}(g_1)(t)\right)^\alpha\,dt
\right)^{1/\alpha}\\
&\le  C(\alpha, n)
\left(\frac{1}{|I|}
\int_{2I} ||g(t)||^{(2-\tilde\epsilon)\alpha}
\,dt\right)^{\frac{1}{(2-\tilde\epsilon)\alpha}}
\sum_{i=1}^n \left(\frac{1}{|I|}
\int_{2I}
||W^{1/2}(t) A_I^{-1}
e_i||^{2}\,dt\right)^ {\frac{1}{2}}
\\
&\le  C_1(\alpha, n)
\left(\frac{1}{|I|}
\int_{2I} ||g(t)||^{(2-\tilde\epsilon)\alpha}
\,dt\right)^{\frac{1}{(2-\tilde\epsilon)\alpha}}
\sum_{i=1}^n \left\langle
 \langle W\rangle^{-1/2}_I
\langle W\rangle_{2I}\langle W\rangle^{-1/2}_I
 e_i, e_i\right\rangle^{\frac 1 2}
\\
&\le  C_3(\alpha, n, Q_{2,E}(W))
\inf_{x\in I} \left(M ||g||^{p_*}\right)^{\frac{1}{p_*}}
(x),
\endalign
$$
where $p_*=(2-\tilde\epsilon)\alpha<2$. We used
the doubling property of $W$:
 $ \langle W\rangle^{-1/2}_I
\langle W\rangle_{2I}\langle W\rangle^{-1/2}_I
\le 2 \frac{Q^2}{\eta^2}$, the inequality which can be
proved in the same way as (5.6).

The last inequality ensures that for any $\tau$,
$\tau\in (0,1)$, using Kolmogorov-type inequalities
we can find a subset $E(\tau, I_0)\subset I_0$,
$|E(\tau,I_0)|\ge
|I_0|-\tau^\alpha|I|\ge (1-3\tau^\alpha)|I_0|$ such that
$$
t\in E(\tau,I_0)\Rightarrow\quad
S^{A_I^{-1}}(g_1)(t)\le
\frac{C_3(\alpha, n, Q_{2,E}(W))}{\tau}
\inf_{x\in I} \left(M ||g||^{p_*}\right)^{\frac{1}{p_*}}
(x).\tag *
$$

Similarly, for   every $\tau$
there exists a set $E(\tau, I_0)$, 
$|E(\tau,I_0)|\ge (1-3\tau^\alpha)|I_0|$ such that
$$
t\in E(\tau,I_0)\Rightarrow\quad
S^{A_I}(f_1)(t)\le
\frac{C(\alpha, n, Q_{2,E}(W))}{\tau}
\inf_{x\in I} \left(M ||f||^{p_*}\right)^{\frac{1}{p_*}}
(x).\tag *
$$
Here we use the same calculations and the fact that
for any $I$ centered at $E$
$$
\langle W\rangle^{1/2}_I
\langle W^{-1}\rangle_{2I}\langle W\rangle^{1/2}_I
\le 2 Q^2.
$$

Now let us work with $f_2, g_2$.   Let $c_I$ be the
center of the square built on $2 I$. Using the
representation
$$
\left(P_+ A_I W^{-1/2} f_2\right)'(z)=
\frac{1}{2\pi i}\int\frac{(A_I W^{-1/2} f_2)(x)}
{(x-z)^2}dx,\quad\text{Im}\ z>0,\tag 5.15
$$
clearly, we obtain for every $t\in I$
$$
\left(\int_{\Gamma_{t,|I|}}
||\left(P_+ A_I W^{-1/2} f_2\right)'(z)||^2\,d A(z)
\right)^{1/2}\le C
\int\frac{\text{Im}\ c_I}{|x-c_I|^2}
||(A_I W^{-1/2} f_2)(x)||\,dx.\tag 5.16
$$
Therefore, using the inverse H\"older inequality (5.12),
we have again
$$\align
&\left(\int_{\Gamma_{t,|I|}}
||\left(P_+ A_I W^{-1/2} f_2\right)'(z)||^2\,d A(z)
\right)^{1/2}\le C
\sum_{i=1}^n\int\frac{\text{Im}\ c_I}{|x-c_I|^2}
|| W^{-1/2}A_I e_i||\,||f_2||\,dx\\
&\le C_1
\sum_{i=1}^n\left(
\int\frac{\text{Im}\ c_I}{|x-c_I|^2}
|| W^{-1/2}A_I e_i||^{2+\epsilon}\,dx
\right)^{\frac{1}{2+\epsilon}}
\left(\int\frac{\text{Im}\ c_I}{|x-c_I|^2}
||f_2||^{2-\tilde\epsilon}\,dx
\right)^{\frac{1}{2-\tilde\epsilon}}\\
&\le C_2
\sum_{i=1}^n\left\langle 
\langle W\rangle^{1/2}_I
\langle W^{-1}\rangle_{c_I}
\langle W\rangle^{1/2}_I
e_i,e_i \right\rangle^{\frac 1 2}
\inf_{x\in I} \left(M|| f||^{2-\tilde\epsilon}(x)
\right)^{\frac{1}{2-\tilde\epsilon}}.
\tag 5.17
\endalign
$$
Here $2+\epsilon$ is close to 2
($\epsilon\le\epsilon_0$). Finally, using Lemma 5.1
we estimate the last sum by a constant:
$$
\left(\int_{\Gamma_{t,|I|}}
||\left(P_+ A_I W^{-1/2} f_2\right)'(z)||^2\,d A(z)
\right)^{1/2}\le C(n, E, Q)
\inf_{x\in I} \left(M|| f||^{2-\tilde\epsilon}(x)
\right)^{\frac{1}{2-\tilde\epsilon}}.
$$
That is
$$
S^{A_I}(f_2)(t)
\le C(n, E, Q)
\inf_{x\in I} \left(M|| f||^{2-\tilde\epsilon}(x)
\right)^{\frac{1}{2-\tilde\epsilon}},
\quad \forall t\in I.\tag *
$$
The same for $S^{A_I^{-1}}(g_2)(t)$.

Combining all (*) inequalities we obtain that with
a suitable $C=C(n, E, W)$
$$
\align
S_{I}(t):=&
\int_{\Gamma_{t, |I|}}
\left|
\langle(P_+ W^{-1/2}f)'(z), (P_+ W^{1/2}g)'(z)
\rangle_{\Bbb C^n}\right|\,d A(z)\\
\le& S^{A_I}(f)(t)S^{A_I^{-1}}(g)(t)\le C^2
\left(M|| f||^{p_*}(t)\right)^{\frac{1}{p_*}}
\left(M|| g||^{p_*}(t)\right)^{\frac{1}{p_*}}
\tag 5.18
\endalign
$$
at least on a quarter of $I_0$. Of course,
$S_{I_0}(t)\le S_{I}(t)$.

In the case $2 I_0\cap E=\emptyset$ we fix an interval
$I$ centered at $E$ such that $I_0\subset I$ and
$\text{dist}(I_0, E)\ge|I|/6$. Let $c_I$ be the center
of the square built on $I$. We can use again a
representation of the form (5.15):
$$
\left(P_+ A_I W^{-1/2} f\right)'(z)=
\frac{1}{2\pi i}\int\frac{(A_I W^{-1/2} f)(x)}
{(x-z)^2}dx,\quad\text{Im}\ z>0,
$$
to obtain an analog of (5.16)
$$
\left(\int_{\Gamma_{t,|I|}}
||\left(P_+ A_I W^{-1/2} f\right)'(z)||^2\,d A(z)
\right)^{1/2}\le C
\int\frac{\text{Im}\ c_I}{|x-c_I|^2}
||(A_I W^{-1/2} f)(x)||\,dx
$$
for all $t\in I_0$.
Continue in this way we get
$$
S^{A_I}(f)(t)
\le C(n, E, Q)
\inf_{x\in I} \left(M|| f||^{p_*}(x)
\right)^{\frac{1}{p_*}},
\quad \forall t\in I_0.
$$
The same for $S^{A_I^{-1}}(g)(t)$. Thus
$$
S_{I}(t)\le C(n, E, Q)^2
\left(M|| f||^{p_*}(t)\right)^{\frac{1}{p_*}}
\left(M|| g||^{p_*}(t)\right)^{\frac{1}{p_*}}
\tag 5.19
$$
everywhere on $I_0$.

Let $B$ be the largest constant in (5.18), (5.19).
We have already chosen $p_*<2$. Now we introduce the
following function $h(t)$:
$$
h(t)=\sup\{h:\ 
S_{h}(t)\le B
\left(M|| f||^{p_*}(t)\right)^{\frac{1}{p_*}}
\left(M|| g||^{p_*}(t)\right)^{\frac{1}{p_*}}\}.
$$
What we proved can be summarized in:
$$\align
\text{if}\  &I_0:\ 2I_0\cap E\not=\emptyset\quad
\text{then}\quad h(t)\ge|I_0|\quad
\text{on a quarter of measure of}\ I_0,
\\
\text{if}\ &I_0:\ 2I_0\cap E=\emptyset\quad
\text{then}\quad h(t)\ge|I_0|\quad \forall t\in I_0.
\endalign
$$
In any case,
$$\align
&\frac 1 4\int_{\Bbb C_+}
\left|
\langle(P_+ W^{-1/2}f)'(z), (P_+ W^{1/2}g)'(z)
\rangle_{\Bbb C^n}\right|\text{Im}\ z\,d A(z)\\
&\le
\int_{\Bbb R}\int_{\Gamma_{t,h(t)}}
\left|
\langle(P_+ W^{-1/2}f)'(z), (P_+ W^{1/2}g)'(z)
\rangle_{\Bbb C^n}\right|\,d A(z)\,dt\\
&\le
 B\int_{\Bbb R}
\left(M|| f||^{p_*}(t)\right)^{\frac{1}{p_*}}
\left(M|| g||^{p_*}(t)\right)^{\frac{1}{p_*}}\,dt\\
&\le
 B\left(\int_{\Bbb R}
\left(M|| f||^{p_*}(t)\right)^{\frac{2}{p_*}}\,dt
\right)^{\frac 1 2}
\left(\int_{\Bbb R}
\left(M|| g||^{p_*}(t)\right)^{\frac{2}{p_*}}\,dt
\right)^{\frac 1 2}\\
&\le
 B C(p_*)\left(\int_{\Bbb R}
|| f||^{2}(t)\,dt \right)^{\frac 1 2}
\left(\int_{\Bbb R}
|| g||^{2}(t)\,dt\right)^{\frac 1 2}
\endalign
$$
because $\frac{2}{p_*}>1$, and we can use the 
Hardy--Littlewood maximal theorem. The theorem is proved.

\enddemo

\enddocument